\documentclass[review=false, screen=true]{acmart}

\usepackage{amsfonts,amsmath,latexsym,amssymb}
\usepackage{algorithm,algpseudocode}
\usepackage{graphicx,xcolor}
\usepackage{hyperref,cleveref}
\usepackage{tikz}
\usetikzlibrary{3d}
\usepackage{pgfplots}
\usepackage{pgfplotstable}
\usepackage{etoolbox}
\usetikzlibrary{patterns}
\usepackage[caption=false,font=footnotesize]{subfig}
\usepackage{url}
\usepackage{ourmacros}

\graphicspath{{./figs/}}

\newcommand{\LineRef}[1]{\hyperref[#1]{Line~\ref{#1}}}
\newcommand{\lineref}[1]{\hyperref[#1]{line~\ref{#1}}}

\pgfplotsset{filter discard warning=false}
\usepgfplotslibrary{external} 
\usetikzlibrary{calc}

\newcommand{\ignore}[1]{}

\newcommand{\ParNMF}{PLANC}

\makeatletter \def\runningfoot{\def\@runningfoot{}} \def\firstfoot{\def\@firstfoot{}} \makeatother 
  
\begin{document}




\title{PLANC: Parallel Low Rank Approximation with Non-negativity Constraints}

\author{Srinivas Eswar}
\authornote{Eswar and Hayashi share first authorship.}
\affiliation{%
  \institution{Georgia Institute of Technology}
  \city{Atlanta}
  \state{GA}
}
\email{seswar3@gatech.edu}

\author{Koby Hayashi}
\authornotemark[1]
\affiliation{%
  \institution{Georgia Institute of Technology}
  \city{Atlanta}
  \state{GA}
}
\email{khayashi9@gatech.edu}

\author{Grey Ballard}
\orcid{0000-0003-1557-8027}
\affiliation{%
  \institution{Wake Forest University}
  \city{Winston-Salem}
  \state{NC}
}
\email{ballard@wfu.edu}

\author{Ramakrishnan Kannan} 
\authornote{This manuscript has been authored by UT-Battelle, LLC under
Contract No. DE-AC05-00OR22725 with the U.S. Department of Energy.
The United States Government retains and the publisher, by accepting
the article for publication, acknowledges that the United States
Government retains a non-exclusive, paid-up, irrevocable, world-wide
license to publish or reproduce the published form of this manuscript,
or allow others to do so, for United States Government purposes.
The Department of Energy will provide public access to these results of
federally sponsored research in accordance with the DOE Public Access
Plan (\url{http://energy.gov/downloads/doe-public-access-plan}).
}
\orcid{0000-0002-5852-4806}
\affiliation{%
  \institution{Oak Ridge National Laboratory}
  \city{Oak Ridge}
  \state{TN}
}  
\email{kannanr@ornl.gov}

\author{Michael A. Matheson}
\authornotemark[2]
\affiliation{%
  \institution{Oak Ridge National Laboratory}
  \city{Oak Ridge}
  \state{TN}
}  
\email{mathesonma@ornl.gov}

\author{Haesun Park}
\affiliation{%
  \institution{Georgia Institute of Technology}
  \city{Atlanta}
  \state{GA}
}
\email{hpark@gatech.edu}


\begin{abstract}
We consider the problem of low-rank approximation of massive dense non-negative  tensor data, for example to discover latent patterns in video and imaging applications.
As the size of data sets grows, single workstations are hitting bottlenecks in both computation time and available memory.
We propose a distributed-memory parallel computing solution to handle massive data sets, loading the input data across the memories of multiple nodes and performing efficient and scalable parallel algorithms to compute the low-rank approximation.
We present a software package called PLANC (Parallel Low Rank Approximation with Non-negativity Constraints), which implements our solution and allows for extension in terms of data (dense or sparse, matrices or tensors of any order), algorithm (e.g., from multiplicative updating techniques to alternating direction method of multipliers), and architecture (we exploit GPUs to accelerate the computation in this work).
We describe our parallel distributions and algorithms, which are careful to avoid unnecessary communication and computation, show how to extend the software to include new algorithms and/or constraints, and report efficiency and scalability results for both synthetic and real-world data sets.
\end{abstract}

%





\maketitle


\section{Introduction}


The CP decomposition, which is also known as CANDECOMP, PARAFAC, and canonical polyadic decomposition, approximates a tensor, or multidimensional array, by a sum of rank-one tensors.
CP is typically used to identify latent factors in data, particularly when the goal is to interpret those hidden patterns, and it is popular within the signal processing, machine learning, and scientific computing communities, among others \cite{SLFHPF2017,AGHKT14,Hackbusch14}.
Enforcing domain-specific constraints on the computed factors can help to identify interpretable components.
We focus in this paper on non-negative dense tensors (when all tensor entries are nonnegative and nearly all of them are positive) and on constraining solutions to have nonnegative entries. Formally, Non-negative CP (NNCP) can be defined as
\SplitN{\label{eqn:nncp}}{
 \min_{\HH^{(i)} \geq 0} & \left\| \T{X} - \sum_{\rank=1}^\Rank \Mn{H}{1}(:,\rank) \circ \cdots \circ \Mn{H}{N}(:,\rank) \right\|^2
}
where $\Mn{H}{1}(:,i) \circ \cdots \circ \Mn{H}{N}(:,i)$ is the outer product of the $i^{th}$ vector from all the $N$ factors that yields a rank-one tensor and $\sum_{\rank=1}^\Rank \Mn{H}{1}(:,\rank) \circ \cdots \circ \Mn{H}{N}(:,\rank)$ results in a sum of $\Rank$ rank-one tensors that approximate the $N$th order input tensor $\T{X}$. 
For example, in imaging and microscopy applications, tensor values often correspond to intensities, and NNCP can be used to cluster and analyze the data in a lower-dimensional space \cite{JC+16}.
In this work, we consider a brain imaging data set that tracks calcium fluorescence within pixels of a mouse's brain over time during a series of experimental trials \cite{KZ+16}.


The kernel computations within standard algorithms for computing NNCP can be formulated as matrix computations, but the complicated layout of tensors in memory prevents the straightforward use of BLAS and LAPACK libraries.
In particular, the matrix formulation of subcomputations involve different views of the tensor data, so no single layout yields a column- or row-major matrix layout for all subcomputations.
Likewise, the parallelization approach for tensor methods is not a straightforward application of parallel matrix computation algorithms.

In developing an efficient parallel algorithm for computing a NNCP of a dense tensor, the key is to parallelize the bottleneck computation known as Matricized-Tensor Times Khatri-Rao Product (MTTKRP), and a different result is required for each mode of the tensor.
The parallelization must load balance the computation, minimize communication across processors, and distribute the results so that the rest of the computation can be performed independently.
In our algorithm, not only do we load balance the computation, but we also compute and store temporary values that can be used across MTTKRPs of different modes using a technique known as dimension trees, significantly reducing the computational cost compared to standard approaches.
Our parallelization strategy also avoids communicating tensor entries and minimizes the communication of factor matrix entries, helping the algorithm to remain computation bound and scalable to high processor counts.

We employ a variety of algorithmic strategies to computing NNCP, from multiplicative updates to alternating direction method of multipliers.
Because the bottleneck computations such as MTTKRP are shared by all update algorithms that compute gradient information, we separate the parallelization strategy for those computations from the (usually local) computations that are unique to each algorithm.
In this paper, we present an open-source software package called Parallel Low-Rank Approximation with Non-negativity Constraints (PLANC) that currently includes six algorithmic options, and we describe how other algorithms can be incorporated into the framework.
PLANC can also be used for non-negative matrix factorization (NMF) with dense or sparse matrices.
The software is available is \url{https://github.com/ramkikannan/planc}.


Some of the material in this paper has appeared in a previous conference paper \cite{BHK18}, including the parallelization strategy described in \cref{sec:paralg} and the dimension tree optimization detailed in \cref{sec:dimtrees}.
We summarize the main contributions of this paper as follows:
\begin{itemize}
		\item presentation and description of the open-source PLANC software package,
		\item utilization of GPUs to alleviate the MTTKRP bottleneck achieving up to $7\times$ speedup over CPU-only execution,
		\item scaling results for runs from 1 to 16384 nodes (250,000+ cores) on the Titan supercomputer,
		\item side by side run-time and convergence results for various update algorithms,
		\item and new results obtained by applying our code to a mouse brain imaging data set.
\end{itemize}

\section{Non-Negative CP and Alternating-Updating Methods}

\begin{algorithm}
\caption{$\CP = \text{NNCP}(\T{X},\Rank)$}
\label{alg:nncp}
\begin{algorithmic}[1]
\Require $\T{X}$ is $I_1\times \cdots \times I_N$ tensor, $\Rank$ is approximation rank
\State \Comment{Initialize data}
\For{$n=2$ to $N$}
	\State Initialize $\Mn{H}{n}$ 
	\State $\Mn{G}{n} = \MnTra{H}{n}\Mn{H}{n}$
\EndFor
\State \Comment{Compute NNCP approximation}
\While{stopping criteria not satisfied}
	\State \Comment{Perform outer iteration}
	\For{$n=1$ to $N$}
	\State \Comment{Compute new factor matrix in $n$th mode}
	\State $\Mn{M}{n} = \text{MTTKRP}(\T{X},\{\Mn{H}{i}\},n)$
		\label{line:MTTKRP}
	\State $\Mn{S}{n} = \Mn{G}{1} \Hada \cdots \Hada \Mn{G}{n-1} \Hada \Mn{G}{n+1} \Hada \cdots \Hada \Mn{G}{N}$
		\label{line:GH}
	\State $\Mn{H}{n} = \text{NNLS-Update}(\Mn{S}{n},\Mn{M}{n})$
		\label{line:NLS}
	\State $\Mn{G}{n} = \MnTra{H}{n}\Mn{H}{n}$
		\label{line:Gn}
	\EndFor
\EndWhile
\Ensure $\T{X} \approx \CP$
\end{algorithmic}
\end{algorithm}

\begin{figure}[htbp]
\begin{center}
\includegraphics[width=4.5in]{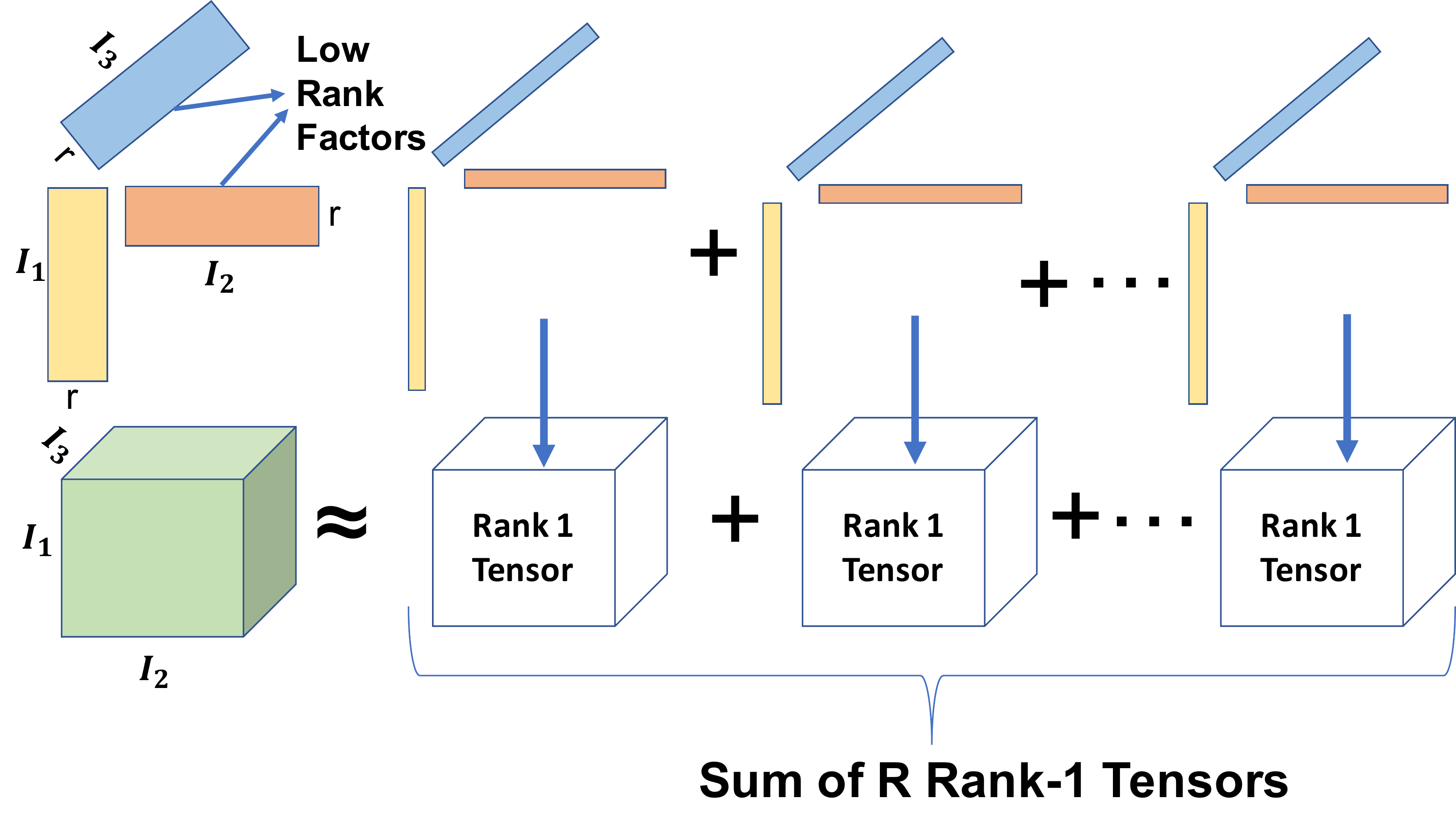}
\caption{CP Decomposition}
\label{fig:cpdecomposition}
\end{center}
\end{figure}

The CP decomposition is a low-rank approximation of a multi-dimensional array, or tensor, which generalizes matrix approximations like the truncated singular value decomposition. 
As in Figure \ref{fig:cpdecomposition}, CP decomposition approximates the given input matrix as sum of $\Rank$ rank-1 tensors.  

\Cref{alg:nncp} shows the pseudocode for an alternating-updating algorithm applied to NNCP \cite{KHP2014}.
\LineRef{line:MTTKRP}, \lineref{line:GH}, and \lineref{line:Gn} compute matrices involved in the gradients of the subproblem objective functions, and \lineref{line:NLS} uses those matrices to update the current factor matrix.

The NNLS-Update in \lineref{line:NLS} can be implemented in different ways.
In a faithful Block Coordinate Descent (BCD) algorithm, the subproblems are solved exactly; in this case, the subproblem is a nonnegative linear least squares problem, which is convex.
We can use the Block Principal Pivoting (BPP) method \cite{KP2011,KHP2014}, which is an active-set-like method, to solve the subproblem exactly.

However, as discussed in \cite{KBP2018} for the matrix case, there are other reasonable alternatives to updating the factor matrix without solving the $N$-block coordinate descent subproblem exactly.
For example, we can more efficiently update individual columns of the factor matrix as is done in the Hierarchical Alternating Least Squares (HALS) method \cite{CP2009}.
In this case, the update rule is 
$$\Mn{H}{n}(:,\rank) \leftarrow \lt[ \Mn{H}{n}(:,\rank) + \Mn{M}{n}(:,\rank) - (\Mn{H}{n} \Mn{S}{n})(:,\rank)  \rt]_+$$
which involves the same matrices $\Mn{M}{n}$ and $\Mn{S}{n}$ as BPP.
Other possible alternating-updating methods include Alternating-Optimization Alternating Direction Method of Multipliers (AO-ADMM) \cite{HSL2015, SBK2017} and Nesterov-based algorithms \cite{LKLHS2017}. The details of each of these  algorithms are presented in Section \ref{sec:nlsalgorithms}. 
The parallel algorithm presented in this paper is generally agnostic to the approach used to solve the nonnegative least squares subproblems, as all these methods are bottlenecked by the subroutine they have in common, the MTTKRP.


\section{Related Work} 
\label{sec:survey}

The formulation of NNCP with least squares error and algorithms for computing it go back as far as~\cite{Paatero97,WW01}, developed in part as a generalization of nonnegative matrix factorization algorithms~\cite{LS99} to tensors.
Sidiropoulos et al.~\cite{SLFHPF2017} provide a more detailed and complete survey that includes basic tensor factorization models with and without constraints, broad coverage of algorithms, and recent driving applications.
The mathematical tensor operations discussed and the notation used in this paper follow Kolda and Bader's survey~\cite{KB2009}. 

Many numerical methods have been developed for the non-negative least squares (NNLS) subproblems which arise in NNCP. Broadly these methods can be divided into projection-based and active-set-based methods.
Projection-based methods are iterative methods which consist of gradient descent and Newton-type algorithms which exploit the fact that the objective function is differentiable and the non-negative projection operator is easy to compute~\cite{lin2007projected,kim2007fast,merritt2005interior,han2009alternating,cichocki2008nonnegative}.
Active-set-like methods explicitly partition the variables into zeros and non-zeros. Once the final partition is known the NNLS problem can be solved via a simpler unconstrained least squares problem~\cite{lawson1974solving,bro1997fast,van2004fast,KP2011}.
We direct the reader to the survey by Kim et al~\cite{KHP2014} for a more in-depth discussion on these methods.

Recently, there has been growing interest in scaling tensor operations to bigger data and more processors in both the data mining/machine learning and the high performance computing communities. 
For sparse tensors, there have been parallelization efforts to compute CP decompositions on shared-memory platforms~\cite{SRSK2015,LCPSV17}, distributed-memory platforms~\cite{KU16,SK16,KU18,MS18} and GPUs~\cite{tang2018gpu,moon2019pl,nisa2019load}, and these approaches can be generalized to constrained problems~\cite{SBK2017}.

Liavas et al.~\cite{LK+17b} extend a parallel algorithm designed for sparse tensors~\cite{SK16} to the 3D dense case.
They use the ``medium-grained'' dense tensor distribution and row-wise factor matrix distribution, which is exactly the same as our distribution strategy (see~\cref{sec:datadist}), and they use a Nesterov-based algorithm to enforce the nonnegativity constraints.
A similar data distribution and parallel algorithm for computing a single dense MTTKRP computation is proposed by Ballard, Knight, and Rouse~\cite{BKR17-TR}. 
Another approach to parallelizing NNCP decomposition of dense tensors is presented by Phan and Cichocki~\cite{PC11}, but they use a dynamic tensor factorization, which performs different, more independent computations across processors. Moon et al.~\cite{moon2019pl} address the data locality optimizations needed during the NNLS phase of the algorithm for both shared memory and GPU systems.

The idea of using dimension trees (discussed in~\cref{sec:dimtrees}) to avoid recomputation within MTTKRPs across modes is introduced in~\cite{PTC13a} for computing the CP decomposition of dense tensors. 
General reuse patterns and mode splitting were present in earlier works on variants of the Tucker Decomposition~\cite{Gra2010,MMVL12}.
It has also been used for sparse CP~\cite{LCPSV17,KU18} and sparse Tucker~\cite{KU16}.


An alternate approach to speeding up CP computations is by reducing the tensor size either via sampling or compression. A large body of work exists for randomized tensor methods~\cite{drineas2011faster,papalexakis2012parcube,wang2015fast,battaglino2018practical} which are recently being extended to the constrained problem~\cite{erichson2018randomized,fu2019block}. The second approach is to first compress the tensor using a different decomposition, like Tucker, and then compute CP on this reduced array. 
This method has been discussed in further detail in~\cite{bro1998improving,tomasi2006comparison}, but it becomes more difficult to impose nonnegative constraints on the overall model.
A separate approach is compute the (constrained) CP decomposition of the entire approximation, rather than only the core tensor, exploiting the structure of the Tucker model to perform the optimization algorithm more efficiently \cite{VDL19}.

%
%


\section{Algorithms} 
\label{sec:algorithm}

\subsection{Parallel NNCP Algorithm}
\label{sec:paralg}

\begin{algorithm}
\caption{$\CP = \text{Par-NNCP}(\T{X},R)$}
\label{alg:Par-NNCP-short}
\begin{algorithmic}[1]
\Require $\T{X}$ is an $I_1\times \cdots \times I_N$ tensor distributed across a $P_1\times \cdots \times P_N$ grid of $P$ processors, so that $\T{X}_{\V{p}}$ is $(I_1/P_1)\times \cdots \times (I_N/P_N)$ and is owned by processor $\V{p}=(p_1,\dots,p_N)$, $\Rank$ is rank of approximation
\For{$n=2$ to $N$}
	\State Initialize $\Mn{H}{n}_{\V{p}}$ of dimensions $(I_n/P)\times \Rank$ 
	\State $\M[\overline]{G} = \text{Local-SYRK}(\Mn{H}{n}_{\V{p}})$
	\State $\Mn{G}{n} = \text{All-Reduce}(\M[\overline]{G},\textsc{All-Procs})$
	\State $\Mn{H}{n}_{p_n} = \text{All-Gather}(\Mn{H}{n}_{\V{p}},\textsc{Proc-Slice}(n,\VE{p}{n}))$
\EndFor
\State \Comment{Compute NNCP approximation}
\While{not converged}
	\label{line:while}
	\For{$n=1$ to $N$}
		\label{line:for}
		\State \Comment{Compute new factor matrix in $n$th mode}
		\State $\M[\overline]{M} = \text{Local-MTTKRP}(\T{X}_{p_1\cdots p_N},\{\Mn{H}{i}_{p_i}\},n)$
			\label{line:locMTTKRP}
		\State $\Mn{M}{n}_{\V{p}} = \text{Reduce-Scatter}(\M[\overline]{M},\textsc{Proc-Slice}(n,\VE{p}{n}))$ 
			\label{line:reduce-scatter}
		\State $\Mn{S}{n} = \Mn{G}{1} \Hada \cdots \Hada \Mn{G}{n-1} \Hada \Mn{G}{n+1} \Hada \cdots \Hada \Mn{G}{N}$
			\label{line:hadamard}
		\State $\Mn{H}{n}_{\V{p}} = \text{NNLS-Update}(\Mn{S}{n},\Mn{M}{n}_{\V{p}})$
			\label{line:locNLS}
		\State \Comment{Organize data for later modes}
		\State $\M[\overline]{G} = {\Mn{H}{n}_{\V{p}}}^\Tra\Mn{H}{n}_{\V{p}}$
			\label{line:locSYRK}
		\State $\Mn{G}{n} = \text{All-Reduce}(\M[\overline]{G},\textsc{All-Procs})$
			\label{line:all-reduce}
		\State $\Mn{H}{n}_{p_n} = \text{All-Gather}(\Mn{H}{n}_{\V{p}},\textsc{Proc-Slice}(n,\VE{p}{n}))$
			\label{line:all-gather}
	\EndFor 
		\label{line:endfor}
\EndWhile
	\label{line:endwhile}
\Ensure $\T{X} \approx \CP$
\Ensure Local matrices: $\Mn{H}{n}_{\V{p}}$ is $(I_n/P)\times \Rank$ and owned by processor $\V{p}=(p_1,\dots,p_N)$, for $1\leq n \leq N$, $\V{\lambda}$ stored redundantly on every processor
\end{algorithmic}
\end{algorithm}

\subsubsection{Algorithm Overview}

The basic sequential algorithm is given in \Cref{alg:nncp}, and the parallel version is given in \Cref{alg:Par-NNCP-short}.
In \Cref{alg:Par-NNCP-short}, we will refer to both the inner iteration, in which one factor matrix is updated (\lineref{line:for} to \lineref{line:endfor}), and the outer iteration, in which all factor matrices are updated (\lineref{line:while} to \lineref{line:endwhile}).
In the parallel algorithm, the processors are organized into a logical multidimensional grid (tensor) with as many modes as the data tensor.
The communication patterns used in the algorithm are MPI collectives: All-Reduce, Reduce-Scatter, and All-Gather.
The processor communicators (across which the collectives are performed) include the set of all processors and the sets of processors within the same processor slice.
Processors within a mode-$n$ slice all have the same $n$th coordinate.

The method of enforcing the nonnegativity constraints of the linear least squares solve or update generally affects only local computation because each row of a factor matrix can be updated independently.
In our algorithm, each processor solves the linear problem or computes the update for its subset of rows (see \lineref{line:locNLS}). 
The most expensive (and most complicated) part of the parallel algorithm is the computation of the MTTKRP, which corresponds to \lineref{line:locMTTKRP}, \lineref{line:reduce-scatter}, and \lineref{line:all-gather}.

The details that are omitted from this presentation of the algorithm include the normalization of each factor matrix after it is computed and the computation of the residual error at the end of an outer iteration.
These two computations do involve both local computation and communication, but their costs are negligible.
We discuss normalization and error computation and give more detailed pseudocode in \Cref{alg:Par-NNCP-long}.


\newcommand{\procdim}{3}
\newcommand{\proc}{\draw[black,shift={(-.5,-.5)}] (0,0) grid (\procdim,\procdim);}
\newcommand{\highlight}{gray!75}
\newcommand{\commhighlight}{gray!25}
\newcommand{\parscale}{.46}
\newcommand{\secfaclabel}{$\Mn{M}{2}$}

\newcommand{\parbasepic}{
\begin{scope}[canvas is yz plane at x=.5,shift={(1.5,-1.5)}]
	\draw[fill=\highlight,shift={(0,1)}] (0,0) rectangle (1,1);
	\draw[fill=\highlight,shift={(-2.5,2-1/3)},xscale=.5] (0,0) rectangle (-1,-1/9);
	\draw[shift={(0,-1.5)},yscale=.5] (0,0) rectangle (1/9,-1);
\end{scope}
\begin{scope}[canvas is zx plane at y=(\procdim-.5),rotate=-90,shift={(-.5,-3)}]
	\draw[fill=\highlight,shift={(0,2.5)}] (0,0) rectangle (1,1);
	\draw[fill=\highlight,yscale=.5] (0,0) rectangle (1/9,-1);
\end{scope}
\begin{scope}[canvas is yx plane at z=.5,yscale=-1,rotate=0]
	\draw[fill=\highlight,shift={(1.5,-.5)}] (0,0) rectangle (1,1);
\end{scope}

\begin{scope}[canvas is yz plane at x=.5,rotate=-90]
	\proc
\end{scope}
\begin{scope}[canvas is yx plane at z=.5,yscale=-1,rotate=0]
	\proc
\end{scope}
\begin{scope}[canvas is zx plane at y=(\procdim-.5),rotate=180]
	\proc
\end{scope}

\begin{scope}[canvas is yz plane at x=.5,shift={(1.5,-.5)}]
	\draw[shift={(-2.5,0)},xscale=.5] (0,-2) grid (-1,1);
	\node[draw=none] at (-3.6,-.5) {\Large $\Mn{H}{1}$};
	\draw[shift={(0,-2.5)},yscale=.5] (-2,0) grid (1,-1);
	\node[draw=none] at (-.5,-3.5) {\Large \secfaclabel};
\end{scope}
\begin{scope}[canvas is zx plane at y=(\procdim-.5),rotate=-90,shift={(1.5,-.5)}]
	\draw[shift={(0,-2.5)},yscale=.5] (-2,0) grid (1,-1);
	\node[draw=none] at (-.5,-3.5) {\Large $\Mn{H}{3}$};
\end{scope}
}
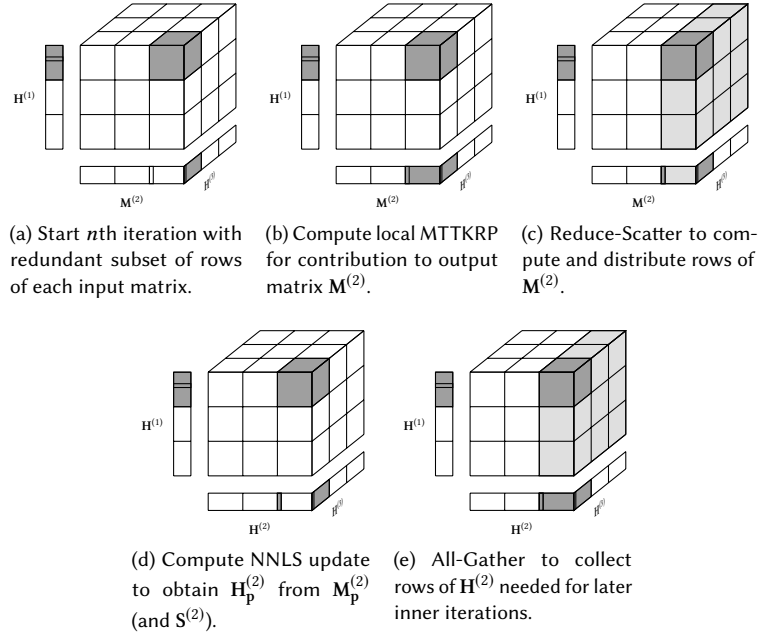
\begin{figure*}[t]
\centering
  \subfloat[Start $n$th iteration with redundant subset of rows of each input matrix. \label{fig:inner_a}]{

\begin{tikzpicture}[x={(-0.5cm,-0.4cm)}, y={(1cm,0cm)}, z={(0cm,1cm)},every node/.append style={transform shape},scale=\parscale]


\begin{scope}[canvas is yz plane at x=.5]
	\draw[fill=\highlight,shift={(-1,.5)},xscale=.5] (0,0) rectangle (-1,-1);
\end{scope}
\begin{scope}[canvas is zx plane at y=(\procdim-.5),rotate=-90]
	\draw[fill=\highlight,yscale=.5,shift={(-.5,-6)}] (0,0) rectangle (1,-1);
\end{scope}

\begin{scope}[canvas is yz plane at x=.5]
\end{scope}

\parbasepic

\end{tikzpicture}}  \quad
  \subfloat[Compute local MTTKRP for contribution to output matrix $\Mn{M}{2}$. \label{fig:inner_b}]{

\begin{tikzpicture}[x={(-0.5cm,-0.4cm)}, y={(1cm,0cm)}, z={(0cm,1cm)},every node/.append style={transform shape},scale=\parscale]


\begin{scope}[canvas is yz plane at x=.5]
	\draw[fill=\highlight,shift={(-1,.5)},xscale=.5] (0,0) rectangle (-1,-1);
\end{scope}
\begin{scope}[canvas is zx plane at y=(\procdim-.5),rotate=-90]
	\draw[fill=\highlight,yscale=.5,shift={(-.5,-6)}] (0,0) rectangle (1,-1);
\end{scope}

\begin{scope}[canvas is yz plane at x=.5]
	\draw[fill=\highlight,shift={(1.5,-3)},yscale=.5] (0,0) rectangle (1,-1);
\end{scope}

\parbasepic

\end{tikzpicture}}  \quad
  \subfloat[Reduce-Scatter to compute and distribute rows of $\Mn{M}{2}$. \label{fig:inner_c}]{

\begin{tikzpicture}[x={(-0.5cm,-0.4cm)}, y={(1cm,0cm)}, z={(0cm,1cm)},every node/.append style={transform shape},scale=\parscale]


\begin{scope}[canvas is yz plane at x=.5]
	\draw[fill=\commhighlight,shift={(1.5,.5)}] (0,0) rectangle (1,-3);
	\draw[fill=\commhighlight,shift={(1.5,-3)},yscale=.5] (0,0) rectangle (1,-1);
	\draw[fill=\highlight,shift={(1.5,-3)},yscale=.5] (0,0) rectangle (1/9,-1);
	\draw[fill=\highlight,shift={(-1,.5)},xscale=.5] (0,0) rectangle (-1,-1);
\end{scope}
\begin{scope}[canvas is zx plane at y=(\procdim-.5),rotate=-90]
	\draw[fill=\commhighlight,shift={(-.5,.5)}] (0,0) rectangle (3,-3);
	\draw[fill=\highlight,yscale=.5,shift={(-.5,-6)}] (0,0) rectangle (1,-1);
\end{scope}
\begin{scope}[canvas is yx plane at z=.5,yscale=-1,rotate=0]
	\draw[fill=\commhighlight,shift={(1.5,-.5)}] (0,0) rectangle (1,3);
\end{scope}

\parbasepic

\end{tikzpicture}} \\
  \subfloat[Compute NNLS update to obtain $\Mn{H}{2}_{\V{p}}$ from $\Mn{M}{2}_{\V{p}}$ (and $\Mn{S}{2}$). \label{fig:inner_d}]{

\renewcommand{\secfaclabel}{$\Mn{H}{2}$}

\begin{tikzpicture}[x={(-0.5cm,-0.4cm)}, y={(1cm,0cm)}, z={(0cm,1cm)},every node/.append style={transform shape},scale=\parscale]


\begin{scope}[canvas is yz plane at x=.5]
	\draw[fill=\highlight,shift={(1.5,-3)},yscale=.5] (0,0) rectangle (1/9,-1);
	\draw[fill=\highlight,shift={(-1,.5)},xscale=.5] (0,0) rectangle (-1,-1);
\end{scope}
\begin{scope}[canvas is zx plane at y=(\procdim-.5),rotate=-90]
	\draw[fill=\highlight,yscale=.5,shift={(-.5,-6)}] (0,0) rectangle (1,-1);
\end{scope}
\begin{scope}[canvas is yx plane at z=.5,yscale=-1,rotate=0]
\end{scope}

\parbasepic

\end{tikzpicture}} \quad
  \subfloat[All-Gather to collect rows of $\Mn{H}{2}$ needed for later inner iterations. \label{fig:inner_e}]{

\renewcommand{\secfaclabel}{$\Mn{H}{2}$}

\begin{tikzpicture}[x={(-0.5cm,-0.4cm)}, y={(1cm,0cm)}, z={(0cm,1cm)},every node/.append style={transform shape},scale=\parscale]


\begin{scope}[canvas is yz plane at x=.5]
	\draw[fill=\commhighlight,shift={(1.5,.5)}] (0,0) rectangle (1,-3);
	\draw[fill=\highlight,shift={(1.5,-3)},yscale=.5] (0,0) rectangle (1,-1);
	\draw[fill=\highlight,shift={(1.5,-3)},yscale=.5] (0,0) rectangle (1/9,-1);
	\draw[fill=\highlight,shift={(-1,.5)},xscale=.5] (0,0) rectangle (-1,-1);
\end{scope}
\begin{scope}[canvas is zx plane at y=(\procdim-.5),rotate=-90]
	\draw[fill=\commhighlight,shift={(-.5,.5)}] (0,0) rectangle (3,-3);
	\draw[fill=\highlight,yscale=.5,shift={(-.5,-6)}] (0,0) rectangle (1,-1);
\end{scope}
\begin{scope}[canvas is yx plane at z=.5,yscale=-1,rotate=0]
	\draw[fill=\commhighlight,shift={(1.5,-.5)}] (0,0) rectangle (1,3);
\end{scope}

\parbasepic

\end{tikzpicture}}
  \caption{Illustration of 2nd inner iteration of Par-NNCP algorithm for 3-way tensor on a $3\times3\times3$ processor grid, showing data distribution, communication, and computation across steps.  Highlighted areas correspond to processor $(1,3,1)$ and its processor slice with which it communicates.  The column normalization and computation of $\Mn{G}{2}$, which involve communication across all processors, is not shown here.}
  \label{fig:inner} 
\end{figure*}

\subsubsection{Data Distribution}
\label{sec:datadist}

Given a logical processor grid of processors $P_1\times \cdots \times P_N$, we distribute the size $I_1 \times \cdots \times I_N$ tensor $\T{X}$ in a block or Cartesian partition.
Each processor owns a local tensor of dimensions $(I_1/P_1)\times \cdots \times (I_N/P_N)$, and only one copy of the tensor is stored.
Locally, the tensor is stored linearly, with entries ordered in a natural mode-descending way that generalizes column-major layout of matrices.
Given a processor $\V{p}=(\VE{p}{1},\dots,\VE{p}{N})$, we denote its local tensor $\T{X}_{\V{p}}$.

Each factor matrix is distributed across processors in a block row partition, so that each processor owns a subset of the rows.
We use the notation $\Mn{H}{n}_{\V{p}}$, which has dimensions $I_n/P\times \Rank$ to denote the local part of the $n$th factor matrix stored on processor $\V{p}$.
However, we also make use a redundant distribution of the factor matrices across processors, because all processors in a mode-$n$ processor slice need access to the same entries of $\Mn{H}{n}$ to perform their computations.
The notation $\Mn{H}{n}_{\VE{p}{n}}$ denotes the $I_n/P_n\times \Rank$ submatrix of $\Mn{H}{n}$ that is redundantly stored on all processors whose $n$th coordinate is $\VE{p}{n}$ (there are $P/P_n$ such processors).

Other matrices involved in the algorithm include $\Mn{M}{n}_{\V{p}}$, which is the result of the MTTKRP computation and has the same distribution scheme as $\Mn{H}{n}_{\V{p}}$, and $\Mn{G}{n}$, which is the $\Rank\times \Rank$ Gram matrix of the factor matrix $\Mn{H}{n}$ and is stored redundantly on all processors.

\subsubsection{Inner Iteration}

The inner iteration is displayed graphically in \Cref{fig:inner} for a 3-way example and an update of the $2$nd factor matrix.
The main idea is that at the start of the $n$th inner iteration (\Cref{fig:inner_a}), all of the data is in place for each processor to perform a local MTTKRP computation, which can be computed using a dimension tree as described in \cref{sec:dimtrees}.
This means that all processors in a slice redundantly own the same rows of the corresponding factor matrix (for all modes except $n$).
After the local MTTKRP is computed (\Cref{fig:inner_b}), each processor has computed a contribution to a subset of the rows of the global MTTKRP $\Mn{M}{n}$, but its contribution must be summed up with the contributions of all other processors in its mode-$n$ slice.
This summation is performed with a Reduce-Scatter collective across the mode-$n$ processor slice that achieves a row-wise partition of the result (in \Cref{fig:inner_c}, the light gray shading corresponds to the rows of $\Mn{M}{2}$ to which processor $(1,3,1)$ contributes and the dark gray shading corresponds to the rows it receives as output).
The output distribution of the Reduce-Scatter is designed so that afterwards, the update of the factor matrix in that mode can be performed row-wise in parallel.
$\Mn{S}{n}$ can be computed locally since the Gram matrices, $\Mn{G}{n}$, are stored redundantly on all processors.
Along with $\Mn{S}{n}$ each processor updates its own rows of the factor matrix given its rows of the MTTKRP result (\Cref{fig:inner_d}).
The remainder of the inner iteration is preparing and distributing the new factor matrix data for future inner iterations, which includes an All-Gather of the newly computed factor matrix $\Mn{H}{n}$ across mode-$n$ processor slices (\Cref{fig:inner_e}) and recomputing $\Mn{G}{n}={\Mn{H}{n}}^\Tra\Mn{H}{n}$.


\Cref{tab:costs} provides the computation, communication, and memory costs of a single outer-iteration, computing $\Mn{S}{n}$ and $\Mn{M}{n}$ for each $n$, which is common to all NNLS algorithms.
We refer the reader to \cite{BHK18} for the detailed analysis of the algorithm and the derivation of these costs.

\begin{table*}
\centering
\begin{tabular}{|ccc|}
\hline
\textbf{Computation} &  \textbf{Communication} & \textbf{Temporary Memory}  \\
\hline
\small
 $O\left(\frac{\Rank}{P} \prod_n I_n + \frac{\Rank^2}{P} \sum_n I_n\right)$ & $\displaystyle O\left(\Rank\sum_n \frac{I_n}{P_n}\right)$ & $\displaystyle O\left(\Rank\left(\prod_n\frac{I_n}{P_n}\right)^{1/2}+\Rank\sum_n \frac{I_n}{P_n}\right)$  \\
\hline
\end{tabular}
\smallskip
\caption{Per-outer-iteration costs in terms of computation (flops), communication (words moved), and memory (words) required to compute $\Mn{S}{n}$ and $\Mn{M}{n}$ for each $n$, assuming the local MTTKRP uses a dimension tree \cite{BHK18}.  These costs do not include the computation (and possibly) communication costs of the particular NNLS algorithm.}
\label{tab:costs}
\end{table*}

\subsection{Dimension Trees}
\label{sec:dimtrees}

\subsubsection{General Approach}

An important optimization of the alternating updating algorithm for NNCP (and unconstrained CP) is to re-use temporary values across inner iterations \cite{PTC13a,KU16-TR,LCPSV17,Kaya17}.
To illustrate the idea, consider a 3-way tensor $\T{X}$ approximated by $\dsquare{\M{U},\M{V},\M{W}}$ and the two MTTKRP computations $\Mn{M}{1}=\underline{\Mz{X}{1}(\M{W}}\Khat\M{V})$ and $\Mn{M}{2}=\underline{\Mz{X}{2}(\M{W}}\Khat\M{U})$ used to update factor matrices $\M{U}$ and $\M{V}$, respectively.
The underlined parts of the expressions correspond to the shared dependence of the outputs on the tensor $\T{X}$ and the third factor matrix $\M{W}$.
Indeed, a temporary quantity, which we refer to as a \emph{partial MTTKRP}, can be computed and re-used across the two MTTKRP expressions.
We refer to the computation that combines the temporary quantity with the other factor matrix to complete the MTTKRP computation as a multi-tensor-times-vector or \emph{multi-TTV}, as it consists of multiple operations that multiply a tensor times a set of vectors, each corresponding to a different mode.

To understand the steps of the partial MTTKRP and multi-TTV operations in more detail, we consider $\T{X}$ to be $I\times J\times K$ and $\M{U}$, $\M{V}$, and $\M{W}$ to have $\Rank$ columns.
Then 
\begin{equation*}
\MnE{M}{1}{i\rank} = \sum_{i,j} \TE{X}{ijk} \ME{V}{j\rank} \ME{W}{k\rank} 
= \sum_{j} \ME{V}{j\rank} \sum_k \TE{X}{ijk} \ME{W}{k\rank} 
= \sum_{j} \ME{V}{j\rank} \TE{T}{ij\rank},
\end{equation*}
where $\T{T}$ is an $I\times J\times \Rank$ tensor that is the result of a partial MTTKRP between tensor $\T{X}$ and the single factor matrix $W$.
Likewise,
\begin{equation*}
\MnE{M}{2}{j\rank} = \sum_{i,k} \TE{X}{ijk} \ME{U}{i\rank} \ME{W}{k\rank} 
= \sum_{i} \ME{U}{i\rank} \sum_k \TE{X}{ijk} \ME{W}{k\rank} 
= \sum_{i} \ME{U}{i\rank} \TE{T}{ij\rank},
\end{equation*}
and we see that the temporary tensor $\T{T}$ can be re-used.
From these expressions, we can also see that computing $\T{T}$ (a partial MTTKRP) corresponds to a matrix-matrix multiplication, and computing each of $\Mn{M}{1}$ and $\Mn{M}{2}$ from $\T{T}$ (a multi-TTV) corresponds to $\Rank$ independent matrix-vector multiplications.
In this case, we compute $\Mn{M}{3}$ using a full MTTKRP.

For a larger number of modes, a more general approach can organize the temporary quantities to be used over a maximal number of MTTKRPs.
The general approach can yield significant benefit, decreasing the computation by a factor of approximately $N/2$ for dense $N$-way tensors.
The idea is introduced in \cite{PTC13a}, but we adopt the terminology and notation of \emph{dimension trees} used for sparse tensors in \cite{KU16-TR,Kaya17}.
In this notation, the root node is labeled $\{1,\dots,N\}$ (we also use the notation $[N]$ for this set) and corresponds to the original tensor, a leaf is labeled $\{n\}$ and corresponds to the $n$th MTTKRP result, and an internal node is labeled by a set of modes $\{i,\dots,j\}$ and corresponds to a temporary tensor whose values contribute to the MTTKRP results of modes $i,\dots,j$.

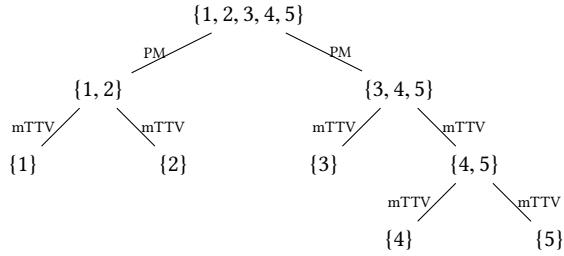
\begin{figure}

\begin{center}
\begin{tikzpicture}

\node (12345) at (4,2) {$\{1,2,3,4,5\}$};
\node (12) at (2,1) {$\{1,2\}$};
\node (345) at (6,1) {$\{3,4,5\}$};
\node (1) at (1,0) {$\{1\}$};
\node (2) at (3,0) {$\{2\}$};
\node (3) at (5,0) {$\{3\}$};
\node (45) at (7,0) {$\{4,5\}$};
\node (4) at (6,-1) {$\{4\}$};
\node (5) at (8,-1) {$\{5\}$};

\scriptsize
\path[draw] (12345) edge [left,align=right] node {PM} (12);
\path[draw] (12345) edge [right,align=left] node {PM} (345);
\path[draw] (12) edge [left,align=right] node {mTTV} (1);
\path[draw] (12) edge [right,align=left] node {mTTV} (2);
\path[draw] (345) edge [left,align=right] node {mTTV} (3);
\path[draw] (345) edge [right,align=left] node {mTTV} (45);
\path[draw] (45) edge [left,align=right] node {mTTV} (4);
\path[draw] (45) edge [right,align=left] node {mTTV} (5);
\normalsize

\end{tikzpicture}
\end{center}
\caption{Dimension tree example for $N=5$. 
The data associated with the root node is the original tensor, the data associated with the leaf nodes are MTTKRP results, and the data associated with internal nodes are temporary tensors.  
Edges labeled with PM correspond to partial MTTKRP computations, and edges labeled with mTTV correspond to multi-TTV computations.}
\label{fig:DT}
\end{figure}

\Cref{fig:DT} illustrates a dimension tree for the case $N=5$.
Various shapes of binary trees are possible \cite{PTC13a,Kaya17}.
For dense tensors, the computational cost is dominated by the root's branches, which correspond to partial MTTKRP computations.
We perform the splitting of modes at the root so that modes are chosen contiguously with the respect to the layout of the tensor entries in memory.
In this way, each partial MTTKRP can be performed via BLAS's GEMM interface without reordering tensor entries in memory.
All other edges in a tree correspond to multi-TTVs and are typically much cheaper.
By organizing the memory layout of temporary quantities, the multi-TTV operations can be performed via a sequence of calls using BLAS's GEMV interface.
By using the BLAS in our implementation, we are able to obtain high performance and on-node parallelism.

\begin{figure}
\subfloat[Partial MTTKRP to compute node $\{3,4,5\}$ from root node $\{1,2,3,4,5\}$, executed via one GEMM call. \label{fig:DT-PM}]{\includegraphics[width=.475\columnwidth]{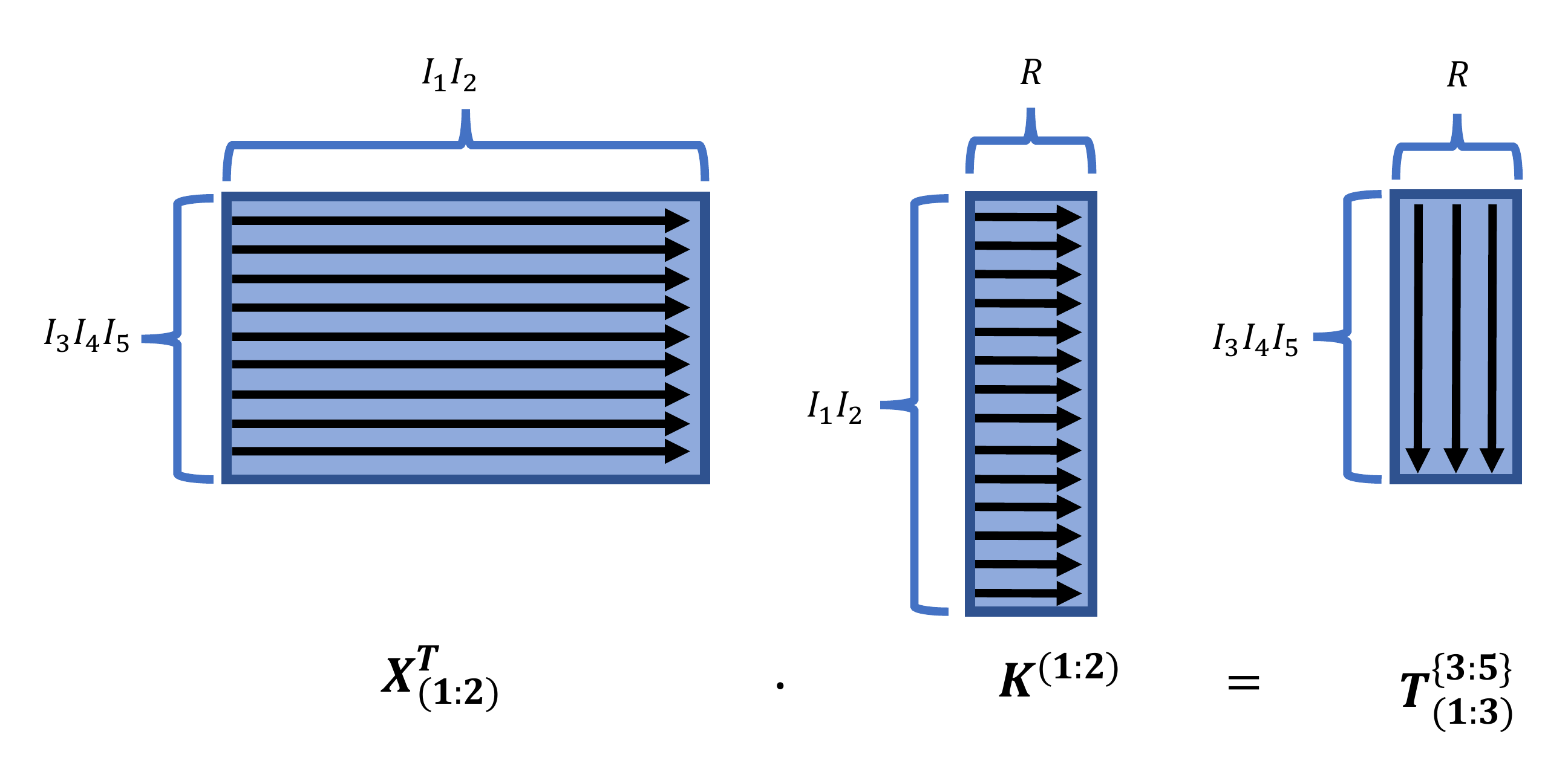}} \hfill
\subfloat[Multi-TTV to compute node $\{3\}$ from node $\{3,4,5\}$, executed via $\Rank$ GEMV calls.  Here $\Mz{T}{1}^{\{3,4,5\}} \lbrack \rank\rbrack$ refers to the $\rank$th contiguous block of size $I_3\times (I_4I_5)$ of $\Mz{T}{1}^{\{3,4,5\}}$. \label{fig:DT-mTTV}]{\includegraphics[width=.475\columnwidth]{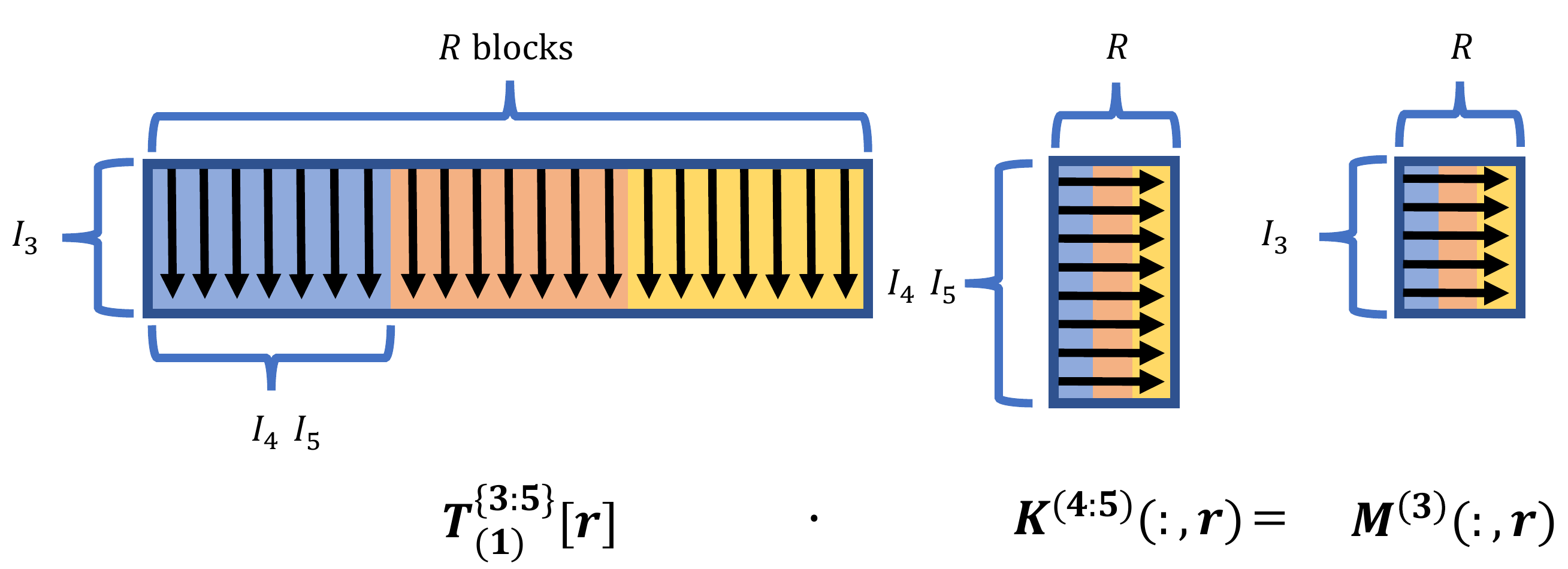}}
\caption{Data layout and dimensions for two example computations in dimension tree shown in \Cref{fig:DT}.  In this notation, $\Mz{X}{1:2}$ is the matricization of input tensor $\T{X}$ with respect to modes 1 through 2, $\Mn{K}{1:2} = \Mn{H}{2} \Khat \Mn{H}{1}$, $\T T^{\{3,4,5\}}$ is the temporary $I_3 \times I_4 \times I_5 \times \Rank$ tensor corresponding to node $\{3,4,5\}$ in the dimension tree, $\Mn{K}{4:5} = \Mn{H}{5} \Khat \Mn{H}{4}$, and $\Mn{M}{3}$ is the MTTKRP result for mode 3.  The arrows represent row- vs column-major ordering in memory.}
\label{fig:DTcartoon}
\end{figure}

\Cref{fig:DTcartoon} shows the data layout and dimensions of a partial MTTKRP and a multi-TTV taken from the example dimension tree in \Cref{fig:DT}.
\Cref{fig:DT-PM} shows a partial MTTKRP between the input tensor $\T{X}$ and the Khatri-Rao product of the factor matrices in modes 1 and 2, which produces a temporary tensor $\T{T}$ corresponding to the $\{3,4,5\}$ node in the dimension tree.
The key to efficiency in this computation is that the matricization of $\T{X}$ that assigns modes 1 through 2 to rows and modes 3 through 5 to columns, which we denote $\Mz{X}{1:2}$, is already column-major in memory.
Thus, we can use the GEMM interface and compute the temporary tensor $\T{T}$ without reordering any tensor entries.
Note that $\T{T}$ is a 4-way tensor in this case, with its last mode of dimension $\Rank$, and the GEMM interface outputs the matrix $\Mz{T}{1:3}$ (where the first three modes are assigned to rows), which is column-major in memory.
\Cref{fig:DT-mTTV} depicts a multi-TTV that computes the result $\Mn{M}{3}$ from $\T{T}$ and the factor matrices in modes 4 and 5.
Here, the tensor $\T{T}$ is matricized with respect to only its first mode (of dimension $I_3$), but this matricization is also column-major in memory.
We choose the ordering of the modes of $\T{T}$ such that each of $\Rank$ contiguous blocks is used to compute one column of the output matrix via a matrix-vector operation with a corresponding column of the Khati-Rao product of the other factor matrices.

No matter how the dimension tree is designed, the computational cost of each partial MTTKRP is $2I\Rank$, where $I$ is the number of tensor entries and $\Rank$ is the rank of the CP decomposition.
This is the same operation count as a full MTTKRP.
The computational cost of a multi-TTV is the number of entries in the temporary tensor, which is the product of a \emph{subset} of the original tensor dimensions multiplied by $\Rank$.
Thus, it is computationally cheaper than the partial MTTKRPs, but it is also memory bandwidth bound.
The other subroutine necessary for implementing the dimension tree approach is the Khatri-Rao product of contiguous sets of factor matrices.
The computational cost of this operation is also typically lower order, but the running time in practice suffers also from being memory bandwidth bound.

\subsubsection{PLANC Implementation}

For a given tensor, it is possible to compute the dimension tree that minimizes overall computation and memory.
However, for most problems, the computation (and actual running time) will be dominated by the choice of split at the root node, and the other split choices will have negligible effect.
The choice of split at the root node has no effect on the computational cost of the two partial MTTKRPs, but it does affect the temporary memory requirement as well as the practical running time, as that split will determine the dimensions of the two GEMM calls.
The three matrix dimensions in the calls are given by the products of the dimensions of the two subsets of modes and the rank of the decomposition.
The amount of additional memory needed is the size of the larger partial MTTKRP result and is $O(I)$ if $\Rank$ is less than the smallest tensor dimension.

To minimize temporary memory and optimize GEMM performance, we seek to split the modes such that the products of each subset of modes are nearly equal.
To respect the memory layout of the tensor, we consider only contiguous subsets of modes, and thus the split depends on only a single parameter $S$, which we refer to as the ``split'' mode, and split the root into nodes $\{1,\dots,S\}$ and $\{S{+}1,\dots,N\}$.
We compute $S$ to be the smallest mode such that the product of the first $S$ modes is greater than the product of the last $N-S$ modes.

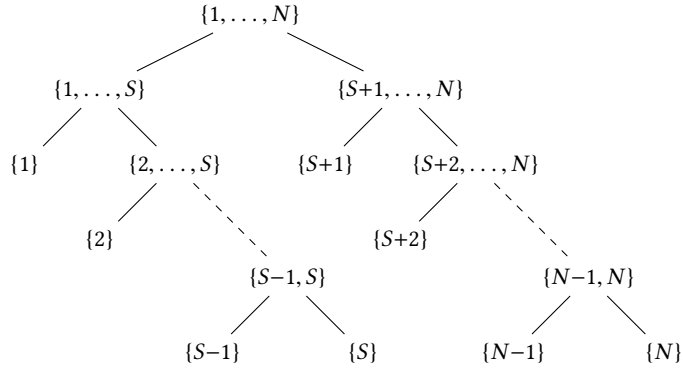
\begin{figure}

\begin{center}
\begin{tikzpicture}

\def\ver{1}
\def\hor{1}

\node (root) at (0,0) {$\{1,\dots,N\}$};
\node (left) at ($ (root) + (-2*\hor,-\ver) $) {$\{1,\dots,S\}$};
\node (right) at ($ (root) + (2*\hor,-\ver) $) {$\{S{+}1,\dots,N\}$};
\node (1) at ($ (left) + (-\hor,-\ver) $) {$\{1\}$};
\node (2S) at ($ (left) + (\hor,-\ver) $) {$\{2,\dots,S\}$};
\node (Sp1) at ($ (right) + (-\hor,-\ver) $) {$\{S{+}1\}$};
\node (Sp2N) at ($ (right) + (\hor,-\ver) $) {$\{S{+2},\dots,N\}$};
\node (Sp2) at ($ (Sp2N) + (-\hor,-\ver) $) {$\{S{+}2\}$};
\node (Nm1N) at ($ (Sp2N) + (1.5*\hor,-1.5*\ver) $) {$\{N{-}1,N\}$};
\node (2) at ($ (2S) + (-\hor,-\ver) $) {$\{2\}$};
\node (Sm1S) at ($ (2S) + (1.5*\hor,-1.5*\ver) $) {$\{S{-}1,S\}$};
\node (Sm1) at ($ (Sm1S) + (-\hor,-\ver) $) {$\{S{-}1\}$};
\node (S) at ($ (Sm1S) + (\hor,-\ver) $) {$\{S\}$};
\node (Nm1) at ($ (Nm1N) + (-\hor,-\ver) $) {$\{N{-}1\}$};
\node (N) at ($ (Nm1N) + (\hor,-\ver) $) {$\{N\}$};

\scriptsize
\path[draw] (root) edge (left);
\path[draw] (root) edge (right);
\path[draw] (left) edge (1);
\path[draw] (left) edge (2S);
\path[draw] (right) edge (Sp1);
\path[draw] (right) edge (Sp2N);
\path[draw] (Sp2N) edge (Sp2);
\path[draw] (Sp2N) edge [dashed] (Nm1N);
\path[draw] (2S) edge (2);
\path[draw] (2S) edge [dashed] (Sm1S);
\path[draw] (Sm1S) edge (Sm1);
\path[draw] (Sm1S) edge (S);
\path[draw] (Nm1N) edge (Nm1);
\path[draw] (Nm1N) edge (N);
\normalsize

\end{tikzpicture}
\end{center}
\caption{Dimension tree used in PLANC for general $N$.  Mode $S$ is the ``split'' mode, chosen so that the product of dimensions in modes $\{1,\dots,S\}$ is approximately equal to that of modes $\{S{+}1,\dots,N\}$.  The splits below the root are not chosen for simplicity.}
\label{fig:ourDT}
\end{figure}

Because the splits within the tree have much less effect on the running time and memory, we structure our tree in order to simplify the software implementation.
That is, we compute the factor matrices in order, from 1 to $N$, and for every internal node of the tree, we split the smallest mode from all other modes.
The structure of the tree we use in PLANC is shown in \Cref{fig:ourDT}, and the pseudocode for its implementation is given by \Cref{alg:ourDT}.
Note that the structure of the main left subtree and the main right subtree are identical, and correspondingly the first half of the pseudocode (for modes 1 to $S$) is nearly identical to the second half (for modes $S{+}1$ to $N$), just with different index ranges.

To explain the pseudocode in more detail, we focus on the first half, or modes 1 through $S$.
The first mode ($n=1$) and the last mode ($n=S$) are special cases because the first mode involves the partial MTTKRP (\lineref{line:n1PM}) and the last mode does not compute an internal node of the tree.
Internal modes ($1<n<S$) involve computing an internal node of the tree and the MTTKRP result for that mode, both of which are computed via multi-TTVs.
We use the notation $\Mn{K}{i:j}$ to represent the reverse Khatri-Rao product of factor matrices $\Mn{M}{i}$ through $\Mn{M}{j}$, which are computed in \lineref{line:n1KRP1}, \lineref{line:n1KRP2}, and \lineref{line:nlSKRP}.
The partial MTTKRP (\lineref{line:n1PM}) is a matrix multiplication between a matricization of the tensor where the first $S$ modes are mapped to rows and a partial Khatri-Rao product; the output is the temporary tensor $\T{T}$, which is computed as a matrix with $\Rank$ columns.
Each matrix involved is either column- or row-major ordered in memory as depicted in \Cref{fig:DT-PM}, for example, where $N=5$.
We use notation $\Mz{T}{1:S}^{\{1:S\}}$ for this output, where the subscript defines the matricization and the superscript labels the temporary tensor corresponding to its node in the dimension tree.
The multi-TTV operations in \lineref{line:n1mTTV}, \lineref{line:nlSmTTV1}, \lineref{line:nlSmTTV2}, and \lineref{line:nSmTTV} are a set of $\Rank$ matrix-vector multiplications.
We use MATLAB-style notation with parentheses to index the $\rank$th column of the Khatri-Rao product matrix and the MTTKRP result matrix.
We use square-bracket notation to index contiguous column blocks of the temporary tensor.
For example, in \lineref{line:nlSmTTV2}, we use $\Mz{T}{1}^{\{n:S\}}[\rank]$ to denote the $\rank$th column block (which comprises $I_{n+1}\cdots I_S$ columns) of the 1st-mode matricization of temporary tensor $\T{T}^{\{n:S\}}$ (which has dimensions $I_n \times \cdots \times I_S \times \Rank$).
This $\rank$th column block is the same as the 1st-mode matricization of the $\rank$th slice of the tensor.
The column blocks are colored distinctly in \Cref{fig:DT-mTTV}, for example, where $\Rank=3$.

We note that for on-node parallelization, we rely on multi-threaded BLAS for the GEMM and GEMV calls, which can be offloaded to a GPU if available.
For the partial Khatri-Rao products, we implement the operation as a row-wise Hadamard product of a set of factor matrix rows, and we use OpenMP parallelization to obtain on-node parallelism.


\begin{algorithm}
\small
\caption{MTTKRP via Dimension Tree}
\label{alg:ourDT}
\begin{algorithmic}[1]
\Require $\T{X}$ is original $N$-way tensor, $\T{T}^{\{i:j\}}$ is temporary tensor of dimension $I_1\times \cdots \times I_{i-1} \times I_{j+1} \times \cdots \times I_N\times R$
\Require $n\in [N]$ is inner iteration mode (evaluated in order), $S\in[N]$ is fixed split mode
\If{$n = 1$}
	\State \label{line:n1KRP1} $\Mn{K}{S{+}1:N} = \Mn{H}{N} \Khat \cdots \Khat \Mn{H}{S+1}$ \hfill \Comment{partial Khatri-Rao product}
	\State \label{line:n1PM} $\Mz{T}{1:S}^{\{1:S\}} = \Mz{X}{1:S} \cdot \Mn{K}{S{+}1:N}$ \hfill \Comment{partial MTTKRP}
	\State \label{line:n1KRP2} $\Mn{K}{1:S{-}1} = \Mn{H}{S-1} \Khat \cdots \Khat \Mn{H}{1}$ \hfill \Comment{partial Khatri-Rao product}
	\State \label{line:n1mTTV} $\Mn{M}{1}(:,r) = \Mz{T}{1}^{\{1:S\}}[r] \cdot \Mn{K}{1:S{-}1}(:,r)$ \; for each $r\in [R]$ \hfill \Comment{multi-TTV for MTTKRP result}
\ElsIf{$n<S$}
	\State \label{line:nlSmTTV1} $\Mz{T}{1:S{-}n{+}1}^{\{n:S\}}(:,r) = \Mz{T}{1}^{\{n{-}1:S\}}[r]^T\cdot \Mn{H}{n-1}(:,r)$\; for each $r\in [R]$ \hfill \Comment{multi-TTV for  internal node tensor}
	\State  \label{line:nlSKRP} $\Mn{K}{n{+}1:S} = \Mz{H}{S} \Khat \cdots \Khat \Mz{H}{n+1}$ \hfill \Comment{partial Khatri-Rao product}
	\State  \label{line:nlSmTTV2} $\Mn{M}{n}(:,r) = \Mz{T}{1}^{\{n:S\}}[r]\cdot \Mn{K}{n{+}1:S}(:,r)$\; for each $r\in [R]$ \hfill \Comment{multi-TTV for MTTKRP result}
\ElsIf{$n=S$}
	\State  \label{line:nSmTTV} $\Mn{M}{S}(:,r) = \Mz{T}{1}^{\{S-1:S\}}[r]\cdot \Mn{H}{S-1}(:,r)$\; for each $r\in [R]$ \hfill \Comment{multi-TTV for MTTKRP result}
\ElsIf{$n=S+1$}
	\State $\Mn{K}{1:S} = \Mn{H}{S} \Khat \cdots \Khat \Mn{H}{1}$ \hfill \Comment{partial Khatri-Rao product}
	\State $\Mz{T}{1:N-S}^{\{S{+}1:N\}} = \Mz{X}{1:S}^T \cdot \Mn{K}{1:S}$ \hfill \Comment{partial MTTKRP}
	\State $\Mn{K}{S{+}2:N} = \Mn{H}{N} \Khat \cdots \Khat \Mn{H}{S+2}$ \hfill \Comment{partial Khatri-Rao product}
	\State $\Mn{M}{S+1}(:,r) = \Mz{T}{1}^{\{S{+}1:N\}}[r] \cdot \Mn{K}{S{+}2:N}(:,r)$ \; for each $r\in [R]$ \hfill \Comment{multi-TTV for MTTKRP result}
\ElsIf{$n<N$}
	\State $\Mz{T}{1:N{-}n{+}1}^{\{n:N\}}(:,r) = \Mz{T}{1}^{\{n{-}1:N\}}[r]^T\cdot \Mn{H}{n-1}(:,r)$\; for each $r\in [R]$ \hfill \Comment{multi-TTV for  internal node tensor}
	\State $\Mn{K}{n{+}1:N} = \Mz{H}{N} \Khat \cdots \Khat \Mz{H}{n+1}$ \hfill \Comment{partial Khatri-Rao product}
	\State $\Mn{M}{n}(:,r) = \Mz{T}{1}^{\{n:N\}}[r]\cdot \Mn{K}{n{+}1:N}(:,r)$\; for each $r\in [R]$ \hfill \Comment{multi-TTV for MTTKRP result}
\Else
	\State $\Mn{M}{N}(:,r) = \Mz{T}{1}^{\{N-1:N\}}[r]\cdot \Mn{H}{N-1}(:,r)$\; for each $r\in [R]$ \hfill \Comment{multi-TTV for MTTKRP result}
\EndIf
\end{algorithmic}
\end{algorithm}

\subsection{Update Algorithms} \label{sec:nlsalgorithms}

In this subsection we consider updating algorithms for the non-negative least squares (NNLS) updates of the factor at each inner iteration of the algorithm (\lineref{line:NLS} of \Cref{alg:nncp}). 
The general problem to be solved in each inner iteration is a constrained least squares problem of the form
\begin{align}
\M{X} \leftarrow \underset{\M{X} \ge \M{0}}{\argmin} \left\lVert \M{A}\M{X} -\M{B}\right\rVert^2_F.
\label{eqn:nls}
\end{align}
All our updating methods (approximately) solve~\Cref{eqn:nls} by first forming $\M{A}^T\M{A}$ and $\M{A}^T\M{B}$, matrices that appear in the gradient of the objective function. 
In the case of updating the factor matrix $\Mn{H}{n}$ we need to solve~\Cref{eqn:nls} with $\M{X} = \MnTra{H}{n}$, $\M{A} = \Mn{K}{n}$, where $\Mn{K}{n}$ is the KRP of factor matrices leaving out the $n^{th}$ factor matrix and $\M{B} = \Mn{X}{n}$, where $\Mn{X}{n}$ is the $n^{th}$ mode matricization of $\T{X}$. 
In this case we have $\M{A}^T\M{A} = \Mn{S}{n}$ and $\M{A}^T\M{B} = {\Mn{M}{n}}^T$, which correspond to the inputs to the NNLS-Update function in \lineref{line:NLS} of \Cref{alg:nncp}.

A nice property of the~\Cref{eqn:nls} is that it can be decoupled along the columns of $\M{X}$ and thus parallelized as in \Cref{alg:Par-NNCP-short}. 
We use the notation $\M{X}_{\V{p}}$ to refer to a subset of the columns of $\M{X}$ owned by processor $\V{p}$, or in the case of \lineref{line:locNLS} of \Cref{alg:Par-NNCP-short}, we use $\Mn{H}{n}_{\V{p}}$ to refer to a subset of the rows of $\Mn{H}{n}=\M{X}^T$.
The gradient for this subset of columns depends on the corresponding columns of $\M{A}^T\M{B} = {\Mn{M}{n}}^T$, denoted by $\Mn{M}{n}_{\V{p}}$, and all of $\M{A}^T\M{A} = \Mn{S}{n}$.


Our framework is capable of supporting any alternating-updating NNCP algorithm \cite{KBP16}. The updating algorithms that fit this framework and are implemented in PLANC are Multiplicative Update~\cite{LS99}, Hierarchical Alternating Least Squares~\cite{Ho2008,CP2009}, Block Principal Pivoting~\cite{KP2011}, Alternating Direction Method of Multipliers~\cite{HSL2016} and Nesterov-type algorithm~\cite{LKLHS2017}. 
We briefly describe the different solvers below. 
Note that the descriptions are for the general form of the NNLS problem as shown in~\Cref{eqn:nls}.

\subsubsection{Multiplicative Update (MU)}
The MU solve is an elementwise operation~\cite{LS99}. The update rule for element $(i,j)$ of $\M{X}$ is
\begin{align}
\M{X}(i,j) \leftarrow \M{X}(i,j)\frac{\M{A}^T\M{B}(i,j)}{\left(\M{A}^T\M{A}\M{X}\right)(i,j)}
\label{eqn:muupd}
\end{align}
While this rule does not solve~\Cref{eqn:nls} to optimality it ensures a reduction in the objective value from the initial value of $\M{X}$. Note that~\Cref{eqn:muupd} breaks down if the denominator becomes zero. In practice a small value is added to the denominator to prevent this situation.

\subsubsection{Hierarchical Alternating Least Squares (HALS)}
HALS updates are performed on individual rows of $\M{X}$~\cite{Ho2008,CP2009}. The update rule for row $i$ can be written in closed form as
\begin{align}
\M{X}(i,:) \leftarrow \left[ \M{X}(i,:) + \frac{(\M{A}^T\M{B})(i,:) - \left(\M{A}^T \M{A}(i,:)\M{X}\right)}{(\M{A}^T \M{A})(i,i)} \right]_+
\label{eqn:halsupd}
\end{align}
where $[\cdot]_+$ is the projection operator onto $\mathbb{R}_+$.
The rows of $\M{X}$ are updated in order so that the latest values are used in every update step. HALS has been observed to produce unbalanced results with either very large or very small values appearing in the factor matrices~\cite{Ho2008,KHP2014}. Normalizing the rows of $\M{X}$ after every update via~\Cref{eqn:halsupd} has been proposed to alleviate this problem~\cite{Ho2008,KHP2014}. 
Within PLANC's parallelization, this step requires explicit communication among processors because the rows of $\M{X}$ (the columns of $\Mn{H}{n}$) are distributed across processors.

\subsubsection{Block Principal Pivoting (BPP)}
BPP is an active-set like method for solving NNLS problems. The main subroutine of BPP is the single right hand side version of~\Cref{eqn:nls},
\begin{align}
\V{x} \leftarrow \underset{\V{x} \ge 0}{\argmin} \left\lVert \M{A}\V{x} -\V{b}\right\rVert^2_2
\label{eqn:snls}
\end{align}
The Karash-Kuhn-Tucker (KKT) optimality conditions for~\Cref{eqn:snls} are specified by $\V{x}$ and $\V{y} = \M{A}^\Tra\M{A}\V{x} - \M{A}^\Tra\V{b}$:  $\V{x}\geq\V{0}$, $\V{y}\geq \V{0}$, and $\V{x} \Hada \V{y} = \V{0}$, where $\Hada$ is the Hadamard product.
The complementary slackness criteria from the KKT conditions forces the support sets, i.e. the non-zero elements, of $\V{x}$ and $\V{y}$ to be disjoint. 
In the optimal solution, the active-set is the set of indices where $x_i = 0$ and the remaining indices are referred to as the passive set. 
Once the active-set is found, we can find the optimal solution to~\Cref{eqn:snls} by  solving an unconstrained least squares problem on the passive set of indices. 
The BPP algorithm attempts to find the active set by greedily swapping indices between the intermediate active and passive sets until it finds a solution that satisfies the KKT conditions. 
The unconstrained least squares is solved using the normal equations. 
Kim and Park discuss the method in greater detail in~\cite{KP2011}. 

\subsubsection{Alternating Direction Method of Multipliers (ADMM)}
In the ADMM solver~\cite{HSL2016} the optimization problem~\Cref{eqn:nls} is reformulated by introducing an auxiliary variable $\M{\hat{X}}$:
\begin{align}
\begin{split}
&\underset{\M{X},\M{\hat{X}}}{\min} \quad \frac 1 2 \left\lVert \M{A}\M{\hat{X}} - \M{B}\right\rVert_F^2 + r(\M{X}) \\
&\textnormal{subject to } \M{X} = \M{\hat{X}}
\end{split}
\label{eqn:admmnls}
\end{align}
where $r(\cdot)$ is the penalty function for nonnegativity. It is 0 if $\M{X} \ge 0$ and $\infty$ otherwise.
The updates for the ADMM algorithm are given by
\begin{align}
\begin{split}
\M{\hat{X}} &\leftarrow \left(\M{A}^\Tra\M{A} + \rho\M{I} \right)^{-1}\left(\M{A}^\Tra\M{B} + \rho\left(\M{X} + \M{U}\right)^\Tra\right) \\
\M{X} &\leftarrow \underset{\M{X}}{\argmin}\ r(\M{X}) + \frac{\rho}{2} \left\lVert \M{X} - \M{\hat{X}} + \M{U}\right\rVert_F^2 \\
\M{U} &\leftarrow \M{U} + \M{X} - \M{\hat{X}},
\end{split}
\label{eqn:admmupd}
\end{align}
where $\M{U}$ is the scaled version of the dual variables corresponding to the equality constraints $\M{X} = \M{\hat{X}}$ and $\rho$ is a step size specified by the user. $\M{U}$ is initialized as a matrix of all zeros.
The advantage of using ADMM is the clever splitting of the non-negativity constraints into updates of two blocks of variables $\M{X}$ and $\M{\hat{X}}$. 
This allows for an unconstrained least squares solve for $\M{\hat{X}}$ and element-wise projections onto $\mathbb{R}_+$ for $\M{X}$. 

We can accelerate this solve by repeating the updates given by~\Cref{eqn:admmupd} more than once. 
One important fact to notice is that the same matrix $\M{A}^\Tra\M{B}$ and matrix inverse $\Br{\M{A}^\Tra\M{A} + \rho \M{I}}^{-1}$ are used for all the updates. 
We can therefore cache $\M{A}^\Tra\M{B}$ and the Cholesky decomposition of  $\Br{\M{A}^\Tra\M{A} + \rho \M{I}}$ to save some computations during subsequent updates. 
We stop updating using the stopping criteria described in~\cite{HSL2016} which is based on $\|\M{X}\|_F$, $\|\M{\hat{X}}\|_F$, and $\|\M{U}\|_F$. 
Computing these norms requires communication because each of these matrices are distributed across processors. 
We also limit the maximum number of acceleration steps to 5.
By default, a good choice for $\rho$ is $\|\M{A}\|_F^2/\Rank$, where $\Rank$ is the number of columns of $\M{A}$ (rank of the CP decomposition)~\cite{HSL2016}.
A comprehensive guide to the ADMM method, convergence properties and selection of optimal $\rho$ can be found in~\cite{boyd2011distributed}.

\subsubsection{Nesterov-type algorithm}
The Nesterov-type algorithm in PLANC was introduced by Liavas et al~\cite{LKLHS2017}. 
Their method solves a modified version of NNLS problem~\Cref{eqn:nls} with the introduction of a proximal term with an auxiliary matrix $\M{X_*}$. 
The proximal term is useful to handle ill-conditioned instances and guarantee strong convexity. 
The objective function tackled is
\begin{align}
f_p(\M{X}) := \frac 1 2 \left\lVert \M{A}\M{X} - \M{B} \right\rVert_F^2 + \frac{\lambda}{2}\left\lVert\M{X} -\M{X_*}\right\rVert_F^2, 
\label{eqn:nesnls}
\end{align}
where $\M{X}$ is constrained to be nonnegative.
The gradient of $f_p$ is given by the expression
\begin{align*}
\nabla f_p(\M{X}) = -(\M{A}^\Tra\M{A}\M{X} - \M{A}^\Tra\M{B}) + \lambda(\M{X}-\M{X_*})
\end{align*}
Updates to $\M{X}$ are performed using the gradient of $f_p$,
\begin{align}
\begin{split}
\nabla f_p(\M{Y}_k) &= \left(\M{A}^\Tra\M{B}-\lambda \M{X_*}\right) + \left(\lambda \M{I} - \M{A}^\Tra\M{A}\right)\M{Y}_k\\
\M{X}_{k+1} &\leftarrow \left[ \M{Y}_k - \alpha \nabla f_p(\M{Y}_k)\right]_+ \\
\M{Y}_{k+1} &\leftarrow \M{X}_{k+1} + \beta_{k+1} \Br{\M{X}_{k+1} - \M{X}_k}
\end{split}
\label{eqn:nesinnerupd}
\end{align}
where $[\cdot]_+$ is the projection operator onto $\mathbb{R}_+$. Notice that we can update $\M{X}$ multiple times reusing $\M{A}^\Tra\M{A}$ and $\M{A}^\Tra\M{B}$. This is the acceleration performed for every inner iteration in \lineref{line:locNLS} of \Cref{alg:Par-NNCP-short}. 
They are repeated until a termination criteria is triggered; different criteria are discussed in~\cite{LKLHS2017}. The termination criteria are bounds checks on the minimum and absolute maximum values of $\M{X}$ and require communication because $\M{X}={\Mn{H}{n}}^\Tra$ is distributed across processors. 
We also limit the total number of inner iterations to 20.

The selection of hyperparameters $\lambda, \alpha$, and $\beta$ depends on the singular values of $\M{A}$  and is necessary for developing a Nesterov-like method for solving~\Cref{eqn:nesnls}. 
The matrix $\M{X_*}$ is generally $\M{X}$ from the previous outer iteration (\lineref{line:while} of \Cref{alg:Par-NNCP-short}).
Details of the selection procedure and different cases can be found in the original paper~\cite{LKLHS2017}. 
 
In addition to the acceleration performed during each NNLS solve,~\Cref{eqn:nesinnerupd}, we can also perform an acceleration step for every outer iteration in the while loop (\lineref{line:while} of \Cref{alg:Par-NNCP-short}). 
In this step all factor matrices are updated using the previous outer iteration values until the objective stops decreasing. The outer acceleration step for iteration $i$ will be,
\begin{align}
\begin{split}
\Mn{H}{1}_{new} &\leftarrow \Mn{H}{1}_i + s_i \Br{\Mn{H}{1}_i - \Mn{H}{1}_{i-1}}\\
\Mn{H}{2}_{new} &\leftarrow \Mn{H}{2}_i + s_i \Br{\Mn{H}{2}_i - \Mn{H}{2}_{i-1}}\\
			&\vdots \\
\Mn{H}{N}_{new} &\leftarrow \Mn{H}{N}_i + s_i \Br{\Mn{H}{N}_i - \Mn{H}{N}_{i-1}}
\end{split}
\label{eqn:nesouterupd}
\end{align}

The results of~\Cref{eqn:nesouterupd} will be accepted as the next iterate only if the overall objective error with the new factor matrices,$\dsquare{\Mn{H}{1}_{new},\dots,\Mn{H}{N}_{new}}$, is lower than that of $\dsquare{\Mn{H}{1}_i,\dots,\Mn{H}{N}_i}$. In order to compute the relative error we need an extra MTTKRP computation per outer acceleration. Typically $s_i = i^{1/N}$ but its value can change as the overall algorithm progresses~\cite{LKLHS2017}.


\section{Software}

\begin{figure}[ht]
\includegraphics[width=\columnwidth]{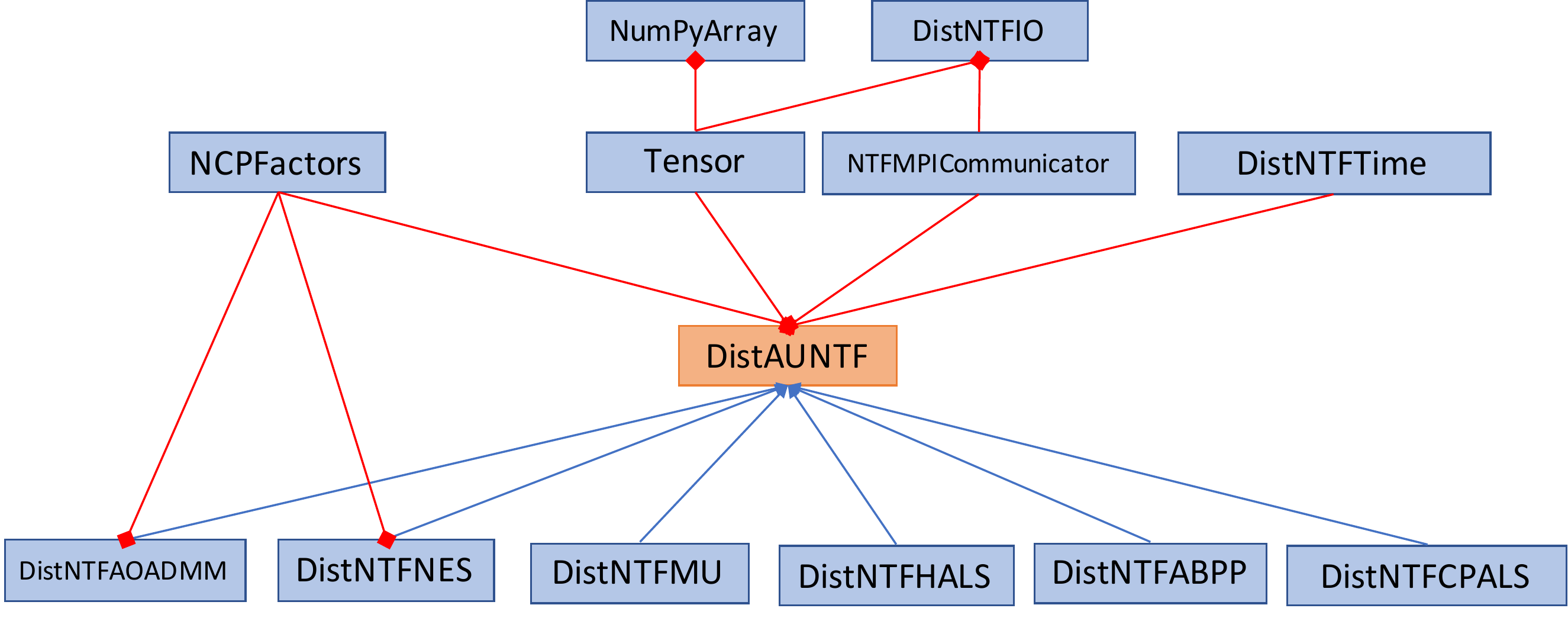}
\caption{\ParNMF\ class diagram. Utility classes are at the top of the diagram, and the algorithm classes at the bottom of the diagram all derive from the abstract class \texttt{DistAUNTF} in orange.  The blue arrows denote an ``is-a'' relationship and the red diamond arrows denote a ``has-a'' relationship.}
\label{fig:classhier}
\end{figure}

The entire \ParNMF\ package has the following modules -- shared memory NMF, shared memory NTF, distributed memory NMF and distributed memory NTF. 
In this section, we give a brief overview of the software package structure for NTF and ways to extend it.

\subsection{Class Organization}
We briefly describe the overall class hierarchy of the \ParNMF\ package as illustrated in \Cref{fig:classhier}. \ParNMF\ offers both shared and distributed memory implementations of NTF and the classes used in each type are distinguished by the prefix \texttt{Dist} in their names (eg. \texttt{DistAUNTF} versus \texttt{AUNTF}). We shall cover the distributed implementation of NTF in this section. Most of the descriptions can be directly applied to the shared memory case as well. 

There are broadly 2 types of classes present. Utility classes are primarily for managing data, setting up the processor grid, and interacting with the user. Algorithm classes perform all the computations needed for NTF and implement the different NNLS solvers.

\subsubsection{Utility Classes}

\paragraph{Data}
The \texttt{Tensor} and \texttt{NCPFactors} classes contain the input tensor $\T{X}$ and the factor matrices $\CP$.
	The \texttt{Tensor} class stores the input tensor as a standard data array. The tensor $\T{X}$ is stored as its mode-1 unfolding $\Mn{X}{1}$ in column major order. Each processor contains  its local part of the tensor (see \Cref{sec:datadist}).
	The \texttt{NCPFactors} class contains all the factor matrices. Each factor matrix is an Armadillo matrix~\cite{sanderson2010}. The matrices are usually column normalized and the column norms are stored in the vector $\V{\lambda}$ which is present as a member of this class (see \Cref{alg:Par-NNCP-long}). The vector $\V{\lambda}$ is replicated in all processors whereas the rows of the factor matrices are distributed across the processor grid (see \Cref{sec:datadist}).
	There is no global view of the entire input tensor or factor matrices and care must be taken to communicate parts of either among the processor grid.

\paragraph{Communication}
The \texttt{NTFMPICommunicator} class creates the MPI processor grid for~\Cref{alg:Par-NNCP-short}. In addition to the communicator for the entire grid, it contains a slice communicator for each mode of the processor grid. The slice communicators are used in the Reduce-Scatter of \lineref{line:reduce-scatter} and All-Gather of \lineref{line:all-gather}  in~\Cref{alg:Par-NNCP-short}.

\paragraph{I/O} 
The \texttt{DistNTFIO} and \texttt{NumPyArray} utility classes are used to read in the input tensor from user-specified files. \texttt{DistNTFIO} also contains methods to generate random tensors and to write out the factor matrices to disk. The \texttt{ParseCommandLine} class contains all the command line options available in \ParNMF. As the name suggests it parses the different combinations of user inputs to instantiate the driver class and run the NTF algorithm. Some example user inputs are the target rank of the decomposition, number of outer iterations, NNLS solver, and regularization parameters.

\subsubsection{Algorithm Classes}

\paragraph{\texttt{DistAUNTF}} 
This is the major workhorse class of the package. It is used to implement~\Cref{alg:Par-NNCP-short}. Some of the important member functions are:
	\begin{itemize}
		 \item \texttt{computeNTF}: This is the outer iteration (\lineref{line:while} in~\Cref{alg:Par-NNCP-short}).
		 \item \texttt{distmttkrp}: Computes the distributed MTTKRP in \lineref{line:locMTTKRP} and \lineref{line:reduce-scatter} in~\Cref{alg:Par-NNCP-short}.
		\item \texttt{gram\_hadamard,update\_global\_gram}: These functions are used to compute the Gram matrix used in the NNLS solvers.
		\item \texttt{computeError}: This function calculates the relative objective error of the factorization as described in~\Cref{sec:error}.
	 \end{itemize}

\paragraph{Derived classes} 
There exist derived classes, such as \texttt{DistNTFANLSBPP}, \texttt{DistNTFMU}, etc., for each of the NNLS solvers described in~\Cref{sec:nlsalgorithms}. There are two main functions which are present in the derived classes which are described below. Auxiliary variables needed to implement certain NNLS solvers like ADMM and Nesterov-type algorithm are also maintained in this class.
	\begin{itemize}
		\item \texttt{update}: This function is the NNLS solve function. It returns the updated factor matrix using the current local MTTKRP result and global Gram matrix (see \Cref{sec:nlsalgorithms}).
		\item \texttt{accelerate}: This implements the outer iteration acceleration (\lineref{line:while} in~\Cref{alg:Par-NNCP-short}). Currently only the Nesterov-type algorithm has an outer  acceleration step.
	\end{itemize}

\subsection{Algorithm Extension}
Extending \ParNMF{} to include different solvers is a simple task and we list the steps to do so below.
\begin{enumerate}
	\item Create a derived class with the newly implemented update function. This is the new NNLS method needed to update the factor matrices.
	\item The constructor for the new class should contain information on whether the algorithm requires an outer acceleration step. If the method requires an outer acceleration step it needs to be implemented in the derived class. 
	\item Update the command line parsing class \texttt{ParseCommandLine} to include additional configuration options for the algorithm.
	\item Include the new algorithm as an option in the utilities and the driver files.
\end{enumerate}

\paragraph{Case Study}
We describe the different steps needed to extend \ParNMF\ to include the Nesterov-type algorithm. 
	\begin{enumerate}
	\item We first create the \texttt{DistNTFNES} class which is derived from \texttt{DistAUNTF}.
	\item We implement the \texttt{update} and \texttt{accelerate} functions in the derived class.
	\begin{itemize}
	\item The Nesterov NNLS updates require the previous iterate values (for the auxiliary term as described in~\Cref{sec:nlsalgorithms}) which may be thought of as a persistent ``state" of the algorithm. We utilize an extra \texttt{NCPFactors} object to hold these variables.
	\item The Nesterov update function needs synchronization in order to terminate its local (iterative) NNLS solve. This involves accessing the communicators found in \texttt{DistNTFMPICommunicator} class for the distributed algorithm. 
	\item Finally, Nesterov-type algorithms generally involve an outer acceleration step which is also implemented in the derived class \texttt{DistNTFNES}.
	\end{itemize}
	\item We then update the \texttt{ParseCommandLine} class to include Nesterov as an algorithm.
	\item We update the driver file \texttt{distntf.cpp} to include the Nesterov algorithm.
\end{enumerate}

\section{Performance Results} \label{sec:experiments}


\newcommand{\datafile}{}
\newcommand{\alg}{}
\newcommand{\numiterations}{30}
\newcommand{\minvalue}{1}
\newcommand{\run}{16}
\newcommand{\nogpu}{0}
\newcommand{\xmax}{10}
\newcommand{\ymin}{0}
\newcommand{\ymax}{10}

\newif\ifnaive
\newif\ifksweep
\newif\iflegend
\newif\ifylabel
\newif\ifwider
\newif\ifliavas
\newif\ifweakthreeD

\newcommand{\setcolors}{
\pgfplotsset{cycle list={
	red, fill=red \\ 
	blue, fill=blue \\ 
	teal, fill=teal \\ 
	red, pattern=crosshatch, pattern color=red \\
	teal, pattern=crosshatch, pattern color=teal \\
	blue, pattern=crosshatch, pattern color=blue \\
}};
}

\newcommand{\breakdownplotoptions}{
	ybar stacked,
	reverse legend,
	bar width=8pt,
	ylabel={Time (s)}, 
	y label style={yshift=-.5cm},
	ymin=0,
	\ifweakthreeD
		symbolic x coords={D1,N1,D8,N8,D27,N27,D64,N64},
	\else
		symbolic x coords={D1,N1,D16,N16,D81,N81,D256,N256},
	\fi	
	xtick=data,
	legend style={at={(0.5,1.3)},anchor=north},
	legend columns=-1,
	reverse legend
}

\newcommand{\breakdownplot}{
\begin{axis}[\breakdownplotoptions]
	\setcolors
	\addplot table[x=alg-p, y expr=(\thisrow{mttkrp}/(\minvalue*\numiterations))] {\datafile};
	\addplot table[x=alg-p, y expr=(\thisrow{krp}/(\minvalue*\numiterations))] {\datafile};
	\addplot table[x=alg-p, y expr=((\thisrow{allgather}+\thisrow{reducescatter})/(\minvalue*\numiterations))] {\datafile};
	\addplot table[x=alg-p, y expr=((\thisrow{nnls}+\thisrow{gram}+\thisrow{allreduce})/(\minvalue*\numiterations))] {\datafile};
	\addplot table[x=alg-p, y expr=((\thisrow{err_compute}+\thisrow{err_communication})/(\minvalue*\numiterations))] {\datafile};
	\legend{MTTKRP,KRP,Factor Comm,NLS,Error};
\end{axis}
}

\newcommand{\labels}{
\node [align=center,text width=3cm] at (1cm, -1.15cm)   {\ifksweep 10 \else 1 \fi};
\node [align=center,text width=3cm] at (2.4cm, -1.15cm)   {\ifksweep 20 \else 16 \fi};
\node [align=center,text width=3cm] at (3.7cm, -1.15cm)   {\ifksweep 30 \else 81 \fi};
\node [align=center,text width=3cm] at (5.1cm, -1.15cm) {\ifksweep 40 \else 256 \fi};
}

\pgfplotsset{
    discard if not/.style 2 args={
        x filter/.append code={
            \edef\tempa{\thisrow{#1}}
            \edef\tempb{#2}
            \ifx\tempa\tempb
            \else
                \def\pgfmathresult{inf}
            \fi
        }
    }
}

\newcommand{\strongscalingplotoptions}{
		log basis y={2},
		log basis x={2},
		\ifliavas
			xlabel=Cores,
		\else
			xlabel=Nodes,
		\fi
		ylabel=Time (s),
		y tick label style={
	        /pgf/number format/.cd,
            	precision=4,
                /tikz/.cd,
		},
	legend style={draw=none, cells={align=left}, nodes={scale=0.7}}
}

\newcommand{\strongscalingplot}{
\begin{loglogaxis}[\strongscalingplotoptions]
	\ifliavas
		\addplot+ [discard if not={alg}{DF},thick,mark options={solid},mark=square*] table [x={p}, y expr=(\thisrow{total}/(\minvalue*\numiterations))] {\datafile};
		\addplot+ [discard if not={alg}{NF},thick,mark options={solid},mark=square*] table [x={p}, y expr=(\thisrow{total}/(\minvalue*\numiterations))] {\datafile};
	\else
		\addplot+ [discard if not={alg}{D},thick,mark options={solid},mark=triangle*] table [x={p}, y expr=(\thisrow{total}/(\minvalue*\numiterations))] {\datafile};
		\addplot+ [discard if not={alg}{N},thick,mark options={solid},mark=square*] table [x={p}, y expr=(\thisrow{total}/(\minvalue*\numiterations))] {\datafile};
	\fi
	\ifliavas
		\addplot+ [discard if not={alg}{L},thick,mark options={solid}] table [x={p}, y expr=(\thisrow{total}/(\minvalue*\numiterations))] {\datafile};
		\legend{DimTree,NoDimTree,NbAO-NTF \cite{LK+17b}}
	\else
		\legend{DimTree,NoDimTree,FlatDimTree,FlatNoDimTree};
	\fi
\end{loglogaxis}
}

\newcommand{\relerrplotoptions}{
	xlabel=Time (s),
	xmin=0,
	xmax=\xmax,
	ymin=\ymin,
	ymax=\ymax,
	\ifylabel
		ylabel=Relative Error,
		ylabel style={yshift=-.3cm},
	\fi
	\ifwider
		width=1.8in, height=2in,
	\else
		width=1.35in, height=2in,
	\fi
	legend style={draw=none, cells={align=left}, nodes={scale=0.7}}
}

\newcommand{\relerrplot}{
\begin{axis}[\relerrplotoptions]
	\addplot+ [unbounded coords=jump,discard if not={algo}{0},discard if not={k}{\run},discard if not={gpu}{1},green,thick,mark options={solid}] table [x={totaltime}, y={relerr}] {\datafile};
	\addplot+ [unbounded coords=jump,discard if not={algo}{1},discard if not={k}{\run},discard if not={gpu}{1},red,thick,mark options={solid},mark=square*] table [x={totaltime}, y={relerr}] {\datafile};
	\addplot+ [unbounded coords=jump,discard if not={algo}{2},discard if not={k}{\run},discard if not={gpu}{1},blue,thick,mark options={solid},mark=triangle*] table [x={totaltime}, y={relerr}] {\datafile};
		\addplot+ [unbounded coords=jump,discard if not={algo}{4},discard if not={k}{\run},discard if not={gpu}{1},magenta,thick,mark options={solid},mark=diamond*] table [x={totaltime}, y={relerr}] {\datafile};
			\addplot+ [unbounded coords=jump,discard if not={algo}{5},discard if not={k}{\run},discard if not={gpu}{1},black,thick,mark options={solid},mark=pentagon*] table [x={totaltime}, y={relerr}] {\datafile};
	\iflegend
		\legend{MU,HALS,ABPP,ADMM,Nesterov}
	\fi
\end{axis}
}

\newcommand{\cpuvsgpuoptions}{
	xlabel=Low Rank,
	xmin=0,
	xmax=\xmax, 
	ymin=0,
	ymax=\ymax,
	\ifylabel
		ylabel=Time(s),
		ylabel style={yshift=-.3cm},
	\fi
	\ifwider
		width=2in, height=2in,
	\else
		width=1.35in, height=2in,
	\fi
	legend style={draw=none, cells={align=left}, nodes={scale=0.7}},
	legend pos={north west}
}

\newcommand{\cpuvsgpuplot}{
\begin{axis}[\cpuvsgpuoptions]
	\addplot+ [discard if not={algo}{0},discard if not={gpu}{\nogpu},discard if not={iter}{\numiterations},green,thick,mark options={solid}] table [x={k}, y={totaltime}] {\datafile};
	\addplot+ [discard if not={algo}{1},discard if not={gpu}{\nogpu},discard if not={iter}{\numiterations},red,thick,mark options={solid},mark=square*] table [x={k}, y={totaltime}] {\datafile};
	\addplot+ [discard if not={algo}{2},discard if not={gpu}{\nogpu},discard if not={iter}{\numiterations},blue,thick,mark options={solid},mark=triangle*] table [x={k}, y={totaltime}] {\datafile};
		\addplot+ [discard if not={algo}{4},discard if not={gpu}{\nogpu},discard if not={iter}{\numiterations},magenta,thick,mark options={solid},mark=diamond*] table [x={k}, y={totaltime}] {\datafile};
			\addplot+ [discard if not={algo}{5},discard if not={gpu}{\nogpu},discard if not={iter}{\numiterations},black,thick,mark options={solid},mark=pentagon*] table [x={k}, y={totaltime}] {\datafile};
	\iflegend
		\legend{MU,HALS,ABPP,ADMM,Nesterov}
	\fi
\end{axis}
}

\newcommand{\lucplot}{
\begin{axis}[\cpuvsgpuoptions]
	\addplot+ [unbounded coords=jump,discard if not={algo}{0},green,thick,mark options={solid}] table [x={k}, y expr=((\thisrow{nnls}/(\minvalue*\numiterations))] {\datafile};
	\addplot+ [unbounded coords=jump,discard if not={algo}{1},red,thick,mark options={solid},mark=square*] table [x={k}, y expr=((\thisrow{nnls}/(\minvalue*\numiterations))] {\datafile};
	\addplot+ [unbounded coords=jump,discard if not={algo}{2},blue,thick,mark options={solid},mark=triangle*] table [x={k}, y expr=((\thisrow{nnls}/(\minvalue*\numiterations))] {\datafile};
		\addplot+ [unbounded coords=jump,discard if not={algo}{4},magenta,thick,mark options={solid},mark=diamond*] table [x={k}, y expr=((\thisrow{nnls}/(\minvalue*\numiterations))] {\datafile};
			\addplot+ [unbounded coords=jump,discard if not={algo}{5},black,thick,mark options={solid},mark=pentagon*] table [x={k}, y expr=((\thisrow{nnls}/(\minvalue*\numiterations))] {\datafile};
	\iflegend
		\legend{MU,HALS,ABPP,ADMM,Nesterov}
	\fi
\end{axis}
}

\subsection{Experimental Setup}

The entire experimentation was performed on Titan, a supercomputer at the Oak Ridge Leadership Computing Facility.  
Titan is a hybrid-architecture Cray XK7 system that contains both advanced 16-core AMD Optero\textsuperscript{TM} central processing units (CPUs) and NVIDIA\textsuperscript{\textregistered} Kepler graphics processing units (GPUs). 
It features 299,008 CPU cores on 18,688 compute nodes, a total system memory of 710 terabytes with 32GB on each node, and Cray's high-performance Gemini network.

We use Armadillo \cite{sanderson2010} for matrix representations and operations.  
In Armadillo, the elements of the dense matrix are stored in column major order.
For dense BLAS and LAPACK operations, we linked Armadillo with the default LAPACK/BLAS wrappers from Cray. 
We use the GNU C++ Compiler (g++ (GCC) 6.3.0) and Cray's MPI library.  
The code can also compile and run on other commodity clusters with entire open source libraries such as OpenBLAS and OpenMPI. 

\subsection{Datasets}

\subsubsection{Mouse Data}
The ``Mouse'' data is a 3D dataset that images a mouse's brain over time and over a sequence of identical trials \cite{KZ+16}.
Each entry of the tensor represents a measure of calcium fluorescence in a particular pixel during a time step of a single trial.
The calcium imaging is performed using an epi-fluorescence macroscope viewing the brain through an artificial crystal skull.
Each image has dimension $1040 \times 1392$, and the minimum number of time steps across 25 trials is 69.
By flattening the pixel dimensions and discarding time steps after 69 for each trial, we obtain a tensor of size $1{,}446{,}680\times 69\times 25$.
Every trial is performed with the same mouse and tracks the same task.
The mouse is presented with visual simulation (starting at frame 3), and after a delay is rewarded with water (starting at frame 25).

\subsubsection{Synthetic}
Our synthetic data sets are constructed from a CP model with an exact low rank with no additional noise.
In this case we can confirm that the residual error of our algorithm with a random start converges to zero.
For the purposes of benchmarking, we run a fixed number of iterations of the NTF algorithms rather than using a convergence check.

\subsection{Performance Breakdown Categories}
The list below gives a brief description of all the categories shown in the breakdown plots and their role in the overall algorithm.

\begin{enumerate}
	\item Gram: the Gram matrix computation includes both the Gram computation of the local factor matrices and the Hadamard product of global Gram matrices for each factor matrix. This computation is performed on each inner iteration but is cheap under the assumption that $\Rank$ is small relative to the tensor dimensions. 
	\item NNLS: the cost of a non-negative least squares update can vary drastically with the algorithm used. The various characteristics that may affect run time for each NNLS algorithm are discussed in \Cref{sec:algdistinct}.
	\item MTTKRP: the (partial) MTTKRP is a purely local computation performed on each node, and can be offloaded to the GPU. Using the dimension tree optimization (\Cref{sec:dimtrees}), we perform 2 partial MTTKRPs for each outer iteration, regardless of the number of modes $N$. Both operations are cast as GEMM calls, where the dimensions are given by the product of the first $S$ mode dimensions ($S$ is the split mode), the product of the last $N-S$ mode dimensions, and the rank $\Rank$.
	\item MultiTTV: the MultiTTVs are purely local computations performed on each node. Each Multi-TTV is cast as a set of $\Rank$ GEMV calls which are typically memory bandwidth bound. 
	\item ReduceScatter: the ReduceScatter collective is used to sum MTTKRP results and distribute portions of the sum appropriately across processors. It is called for each inner iteration.
	\item AllGather: the AllGather collective is used to collect the updated factor matrices to each processor in the slice corresponding to the mode being updated. It is called after each inner iteration.
	\item AllReduce: the AllReduce is used to compute the Gram matrices and for computing norms and other quantities required for stopping criteria of some algorithms. 
\end{enumerate}

\subsection{Updating Algorithm Distinctions}
\label{sec:algdistinct}

\Cref{tab:NNLS_Characteristics} highlights the distinct aspects of each updating algorithm that can affect performance.
The rows of \Cref{tab:NNLS_Characteristics} denote the different local update algorithms implemented in PLANC. 
The algorithms names and acronyms in order from top to bottom in \Cref{tab:NNLS_Characteristics} are as follows: Unconstrained CP (UCP), Multiplicative Update (MU), Hierarchical Least Squares (HALS), Block Principle Pivoting (BPP), Alternating Direction Method of Multipliers (ADMM), and Nesterov-type algorithm (NES).
The aspects of each algorithm that are displayed in \Cref{tab:NNLS_Characteristics} are as follows:
\begin{itemize}
	\item Communication: a check mark and description in this column indicates that the local update algorithm requires some amount of communication. For example, the HALS algorithm requires the communication of the updated column norms. Additional communication requirements can affect performance by incurring additional latency and bandwidth costs. These penalties become significant when the number of processors is high.
	\item Extra MTTKRPs: the NES algorithm has an acceleration step which requires an additional MTTKRP to be performed. This can potentially increase the run time if the acceleration step does not decrease the objective function. 
	Experimentally, on both real and synthetic data sets, we observe that NES run time is significantly increased by the additional MTTKRP computations.
	\item Iteration: this column indicates the iterative nature of the local update algorithm. Note that all of the algorithms we present here are iterative in terms of the outer iteration. UCP, MU, and HALS all have closed-form formulas for the inner iteration, meaning the number of flops can be explicitly computed as a function of the problem size. The rest of the algorithms have flop and communication requirements dependent on the number of iterations it takes the algorithm to converge for a particular local update.
	\item Tuning: many optimization algorithms require some tunable input parameters which can impact performance. For example,   setting a step size is a frequent requirement for gradient based optimization algorithms when an exact line search is too computationally expensive.
\end{itemize}

\begin{table*}
\begin{center}
\begin{tabular}{|c | c | c | c | c|}
\hline
Alg & Communication & Extra MTTKRPs & Iterative & Tuning \\
\hline
\hline
UCP & $\times$ & $\times$ & $\times$  & $\times$ \\
\hline
MU & $\times$  & $\times$  & $\times$ & $\times$\\
\hline
HALS & $\checkmark$ Column Norms& $\times$ & $\times$ & $\times$ \\
\hline
BPP & $\times$ & $\times$ & $\checkmark$ & $\times$  \\
\hline
ADMM & $\checkmark$ Stopping Criteria & $\times$  &$\checkmark$  & $\checkmark$ Step Size \\
\hline
NES & $\checkmark$ Stopping Criteria & $\checkmark$ & $\checkmark$ & $\checkmark$ See \cref{sec:nlsalgorithms} \\
\hline
\end{tabular}
\end{center}
\caption{Characteristics of the various update algorithms that can potentially affect performance. The columns are as follows: 1) if the local update requires communication, 2) if the update requires additional MTTKRP computations, 3) if the local update itself is iterative, 4) if the algorithm's performance are significantly impacted by parameter tuning.  A $\checkmark$ corresponds to the algorithm having the characteristic, and a $\times$ means it does not.}
\label{tab:NNLS_Characteristics}
\end{table*}

\subsection{Microbenchmarks}


\subsubsection{Per-Iteration Timing Comparison across Algorithms}

\begin{figure}
\subfloat[Synthetic data: $384 {\times} 384 {\times} 384 {\times} 384$ tensor on 81 nodes arranged in $3{\times}3{\times}3{\times}3$ grid \label{fig:luclr81}]{
\includegraphics[width=0.45\textwidth, height=2in]{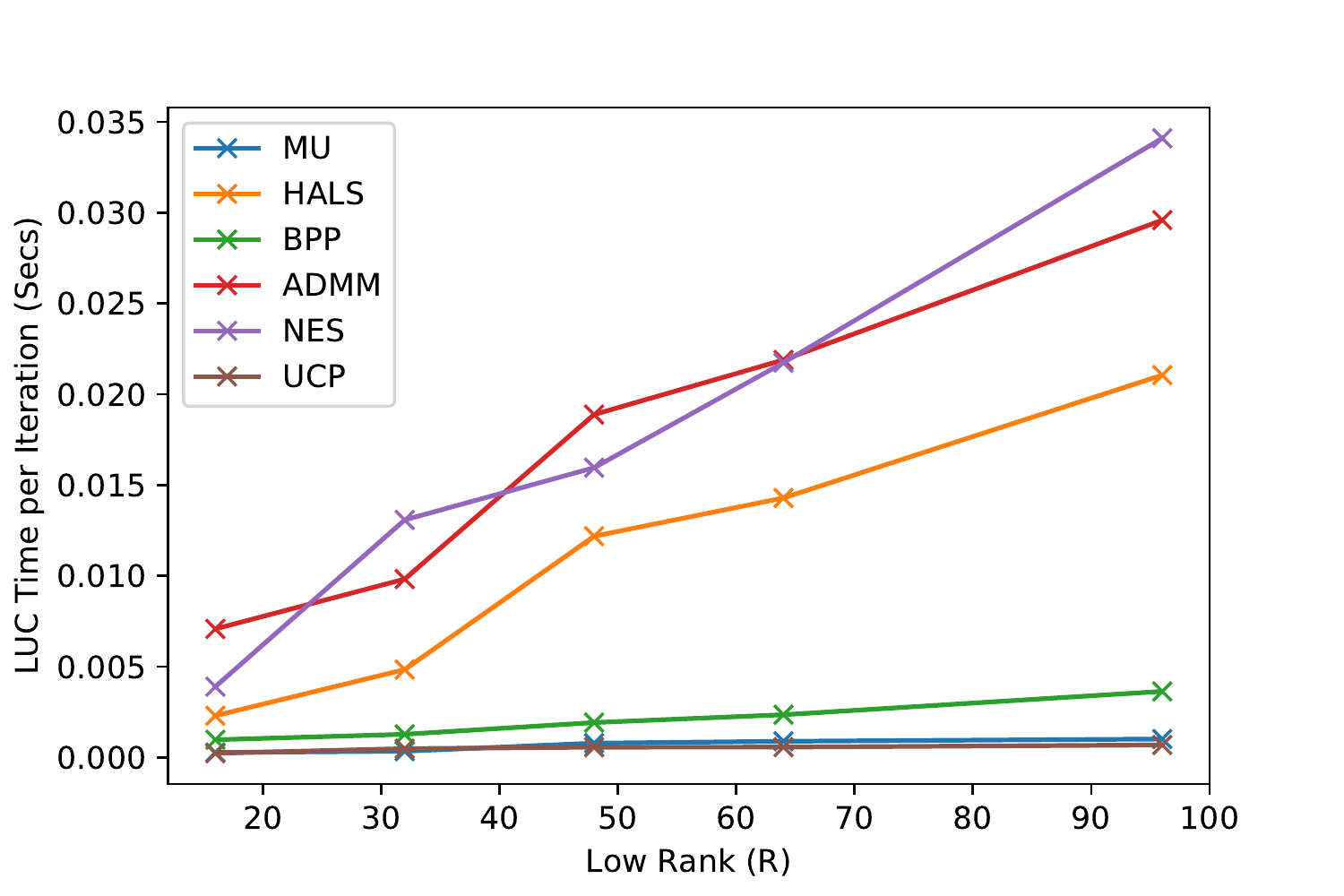}
} \quad \quad
\subfloat[Mouse data: $1447680\times69\times25$ tensor on 64 nodes arranged in $64\times1\times1$ grid \label{fig:lucrw}]{
\includegraphics[width=0.45\textwidth, height=2in]{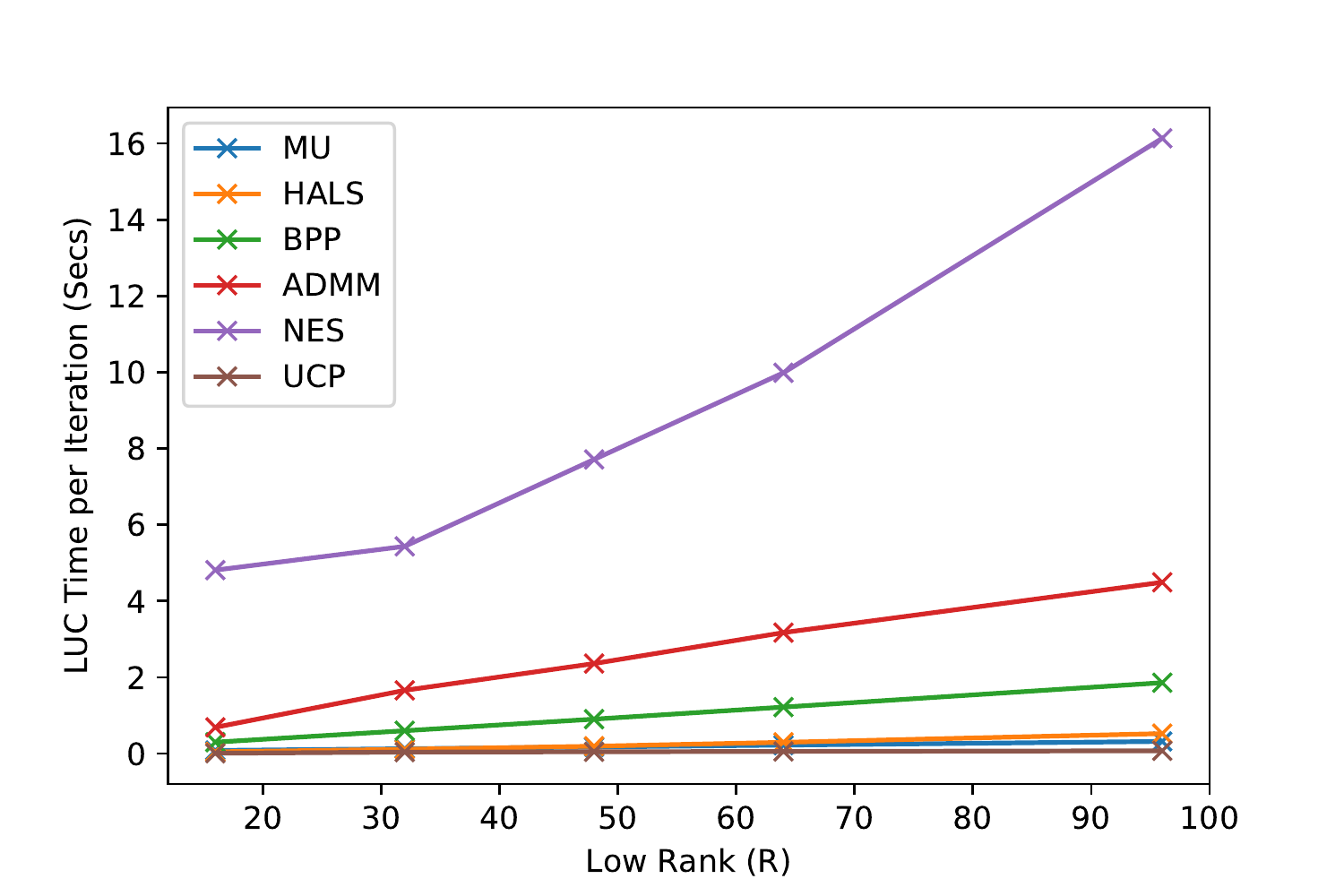}
}
\caption{Per Iteration Local Update Computation (LUC) comparison of NNLS algorithms on 4D synthetic and and 3D Mouse tensors}
\label{fig:synluc}
\end{figure}

\Cref{fig:synluc} shows the ``local update computation" time taken by the different updating algorithms for various low-rank values on a synthetic data set and the Mouse dataset.
The synthetic tensor (\Cref{fig:luclr81}) involves about 20 LUCs per processor per iteration whereas the Mouse data (\Cref{fig:lucrw}) has about 20,000 updates per processor per iteration, which accounts for the difference in the scales of the time seen in the figures.
MU and CP are the cheapest algorithms with NES being the most expensive. HALS, ADMM and NES algorithms all communicate in their update steps and this significantly affects their runtimes, see \Cref{fig:strsca4db}. NES has the most expensive inner iteration involving an eigendecomposition of the Gram matrix and up to 20 iterations of the NNLS updater. ADMM has the second most expensive inner iteration with up to 5 iterations of the acceleration step. HALS on the other hand doesn't have a very expensive inner iteration but needs a synchronization to normalize every updated column of the factor matrix before proceeding to the next column, causing a slowdown.

\subsubsection{Comparison across Processor Grids}

\begin{figure}
\includegraphics[width=\textwidth]{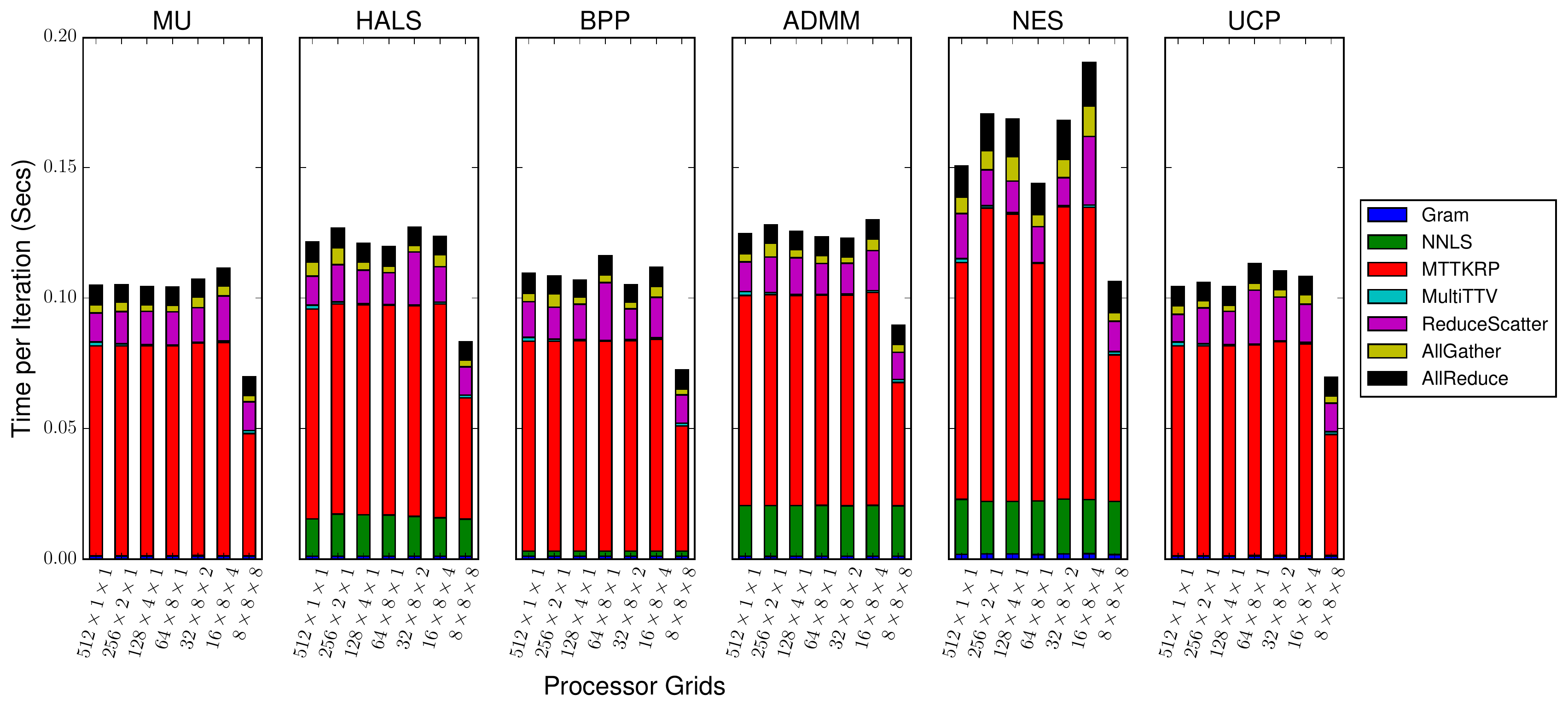}
\caption{Processor grid sweep of a $512 \times 512 \times 512$ synthetic 3D low rank tensor on 512 nodes with low rank 96.}
\label{fig:confsweep3D}
\end{figure}

\Cref{fig:confsweep3D} gives a processor grid comparison for a 3-D cubical tensor of size 512.
The distributed MTTKRP time dominates the overall run time and we observe that an even processor distribution results in the best achieved performance for all update algorithms.
This difference in run times is partially accounted for by GEMM performance due to the different shapes of the matrices involved.
Besides choosing an even processor grid we see that configurations 1 through 6 have quite stable run times with the exception of the NES algorithm.
The variation in the NES run times can be attributed to the variable number of MTTKRPs needed for the acceleration steps.

\subsubsection{Comparison between CPU and GPU matrix multiplication offloading}

\begin{figure}
\centering
\includegraphics[width=0.65\textwidth, height=2.25in]{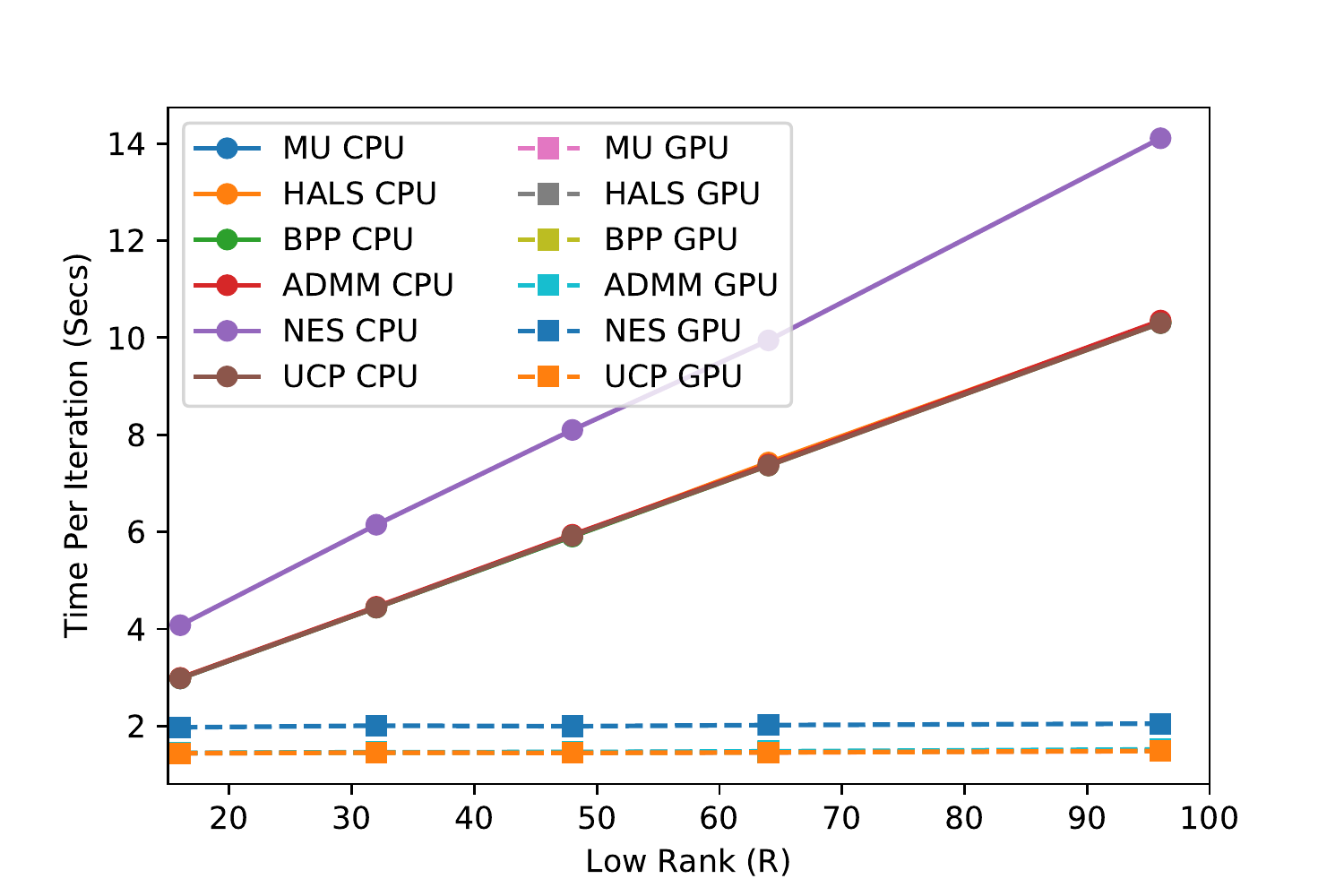}
\caption{Timing comparison of CPU and GPU offloading on 4D Synthetic Low Rank Tensor of size $384{\times}384{\times}384{\times}384$ on $3{\times}3{\times}3{\times}3$ processor grid with varying ranks}
\label{fig:cpuvsgpulowrank}
\end{figure}

\Cref{fig:cpuvsgpulowrank} shows comparison in run times between performing partial MTTKRPs on the CPU versus offloading to the GPU as rank increases.
As expected the CPU run times increase linearly with $\Rank$ as the operation count for CP-ALS is dominated by the MTTKRP which is linear in $\Rank$.
In the case of the GPU execution, the tested sizes of $\Rank$ are never large enough to saturate the GPU, yielding flat run times even as the rank increases.
The NES algorithm takes additional time for both the CPU and GPU executions due to the additional MTTKRPs.
In this case, for the chosen tensor size and rank, it is always beneficial to offload the GEMM calls to the GPU, and the maximum achieved speedup with GPU offloading is about $7\times$.
However, we have observed in other experiments that NVBLAS can make the incorrect decision to offload the computation to the GPU when it is faster to perform it on the CPU.

\subsection{Convergence comparison across algorithms}

\begin{figure}
\centering
\subfloat[Low rank $\Rank=64$ \label{fig:accuracylr81k64}]{
\includegraphics[width=0.45\textwidth, height=2in]{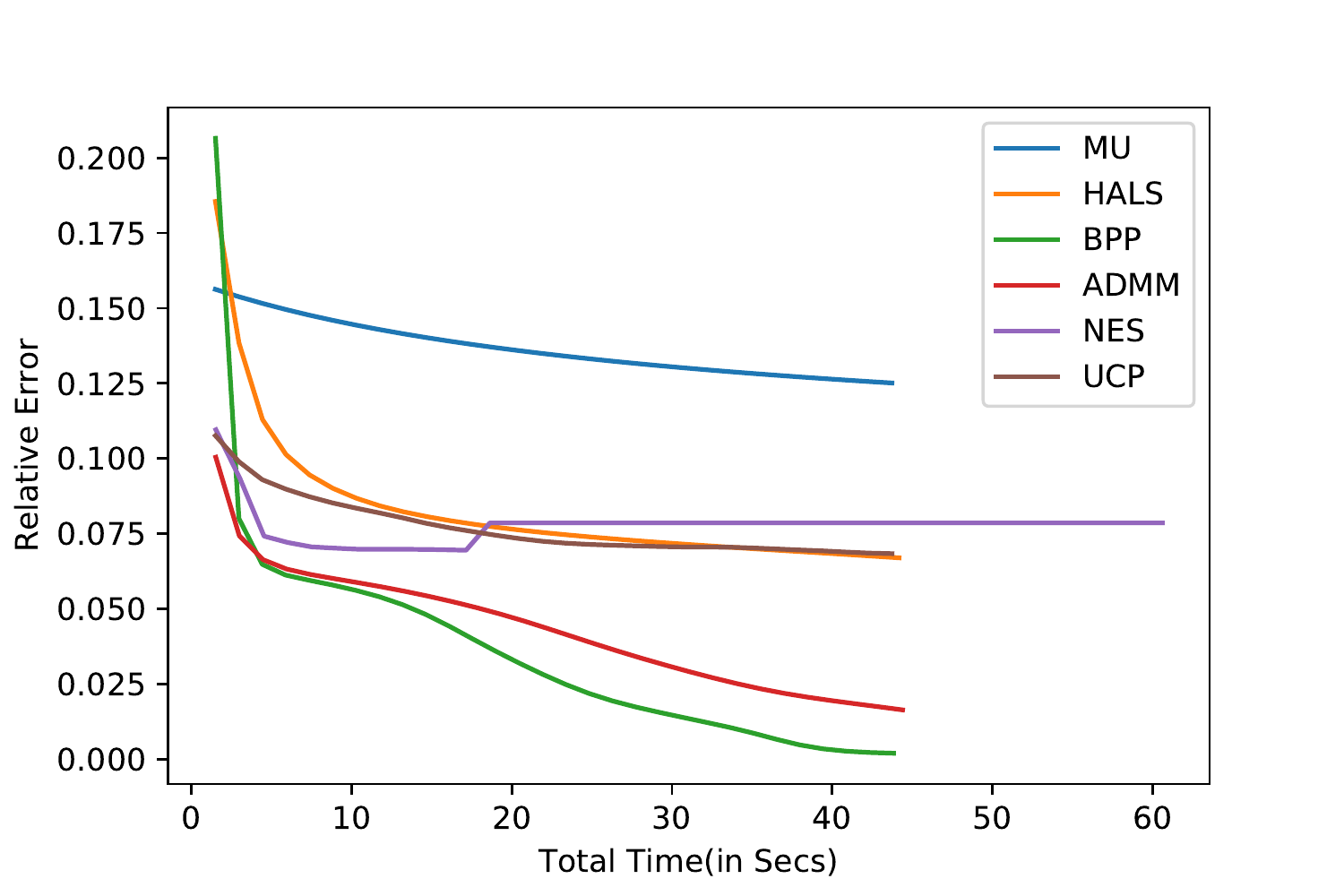}
}
\subfloat[Low rank $\Rank=96$ \label{fig:accuracylr81k96}]{
\includegraphics[width=0.45\textwidth, height=2in]{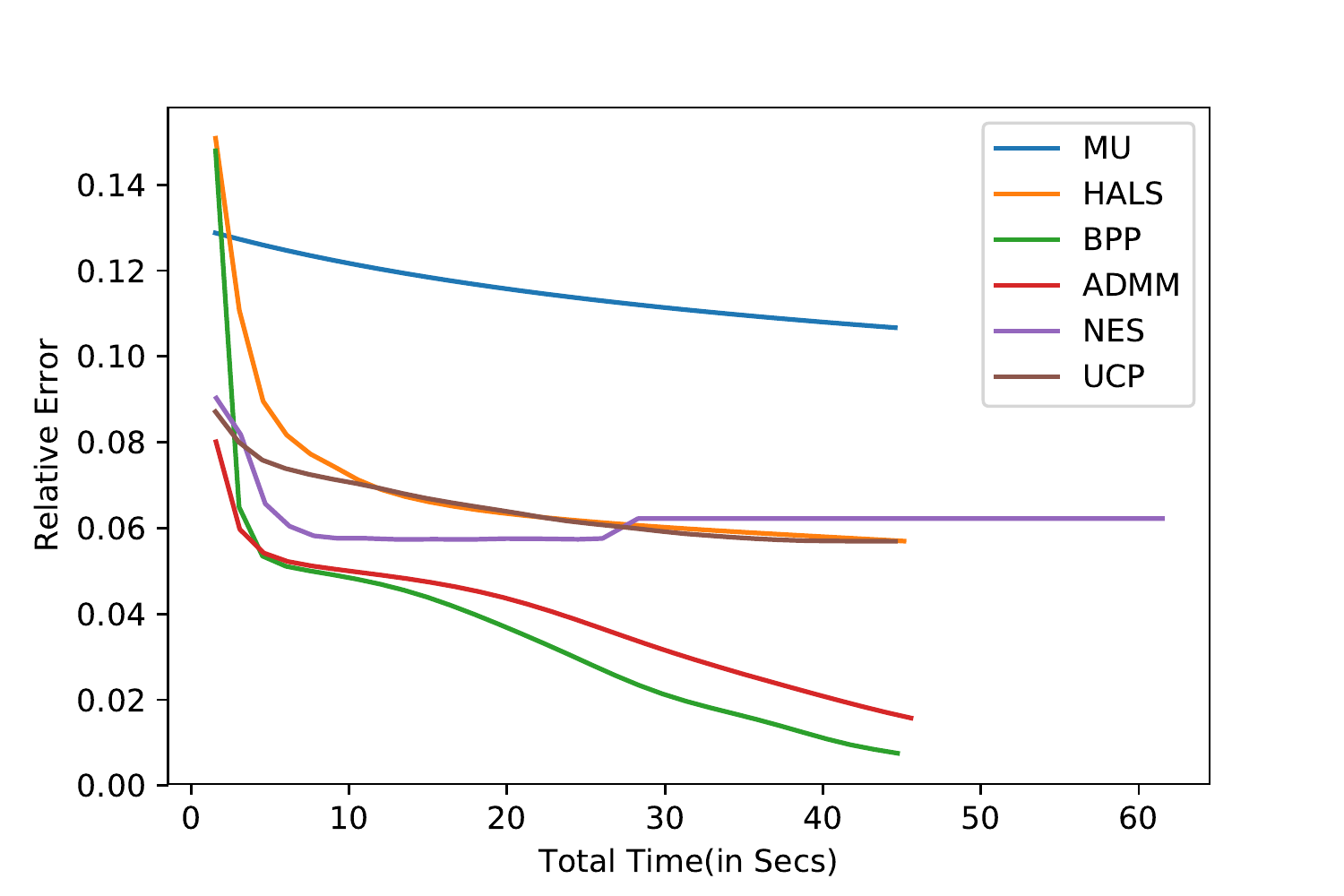}
}
\caption{Relative error over time comparison of updating algorithms on 4D Synthetic Low Rank Tensor of size $384{\times}384{\times}384{\times}384$ on $3{\times}3{\times}3{\times}3$ processor grid}
\label{fig:synconv}
\end{figure}

\begin{figure}
\centering
	\subfloat[Low rank $\Rank=16$ \label{fig:accuracyrwk16}]{
		\includegraphics[width=0.45\textwidth, height=2in]{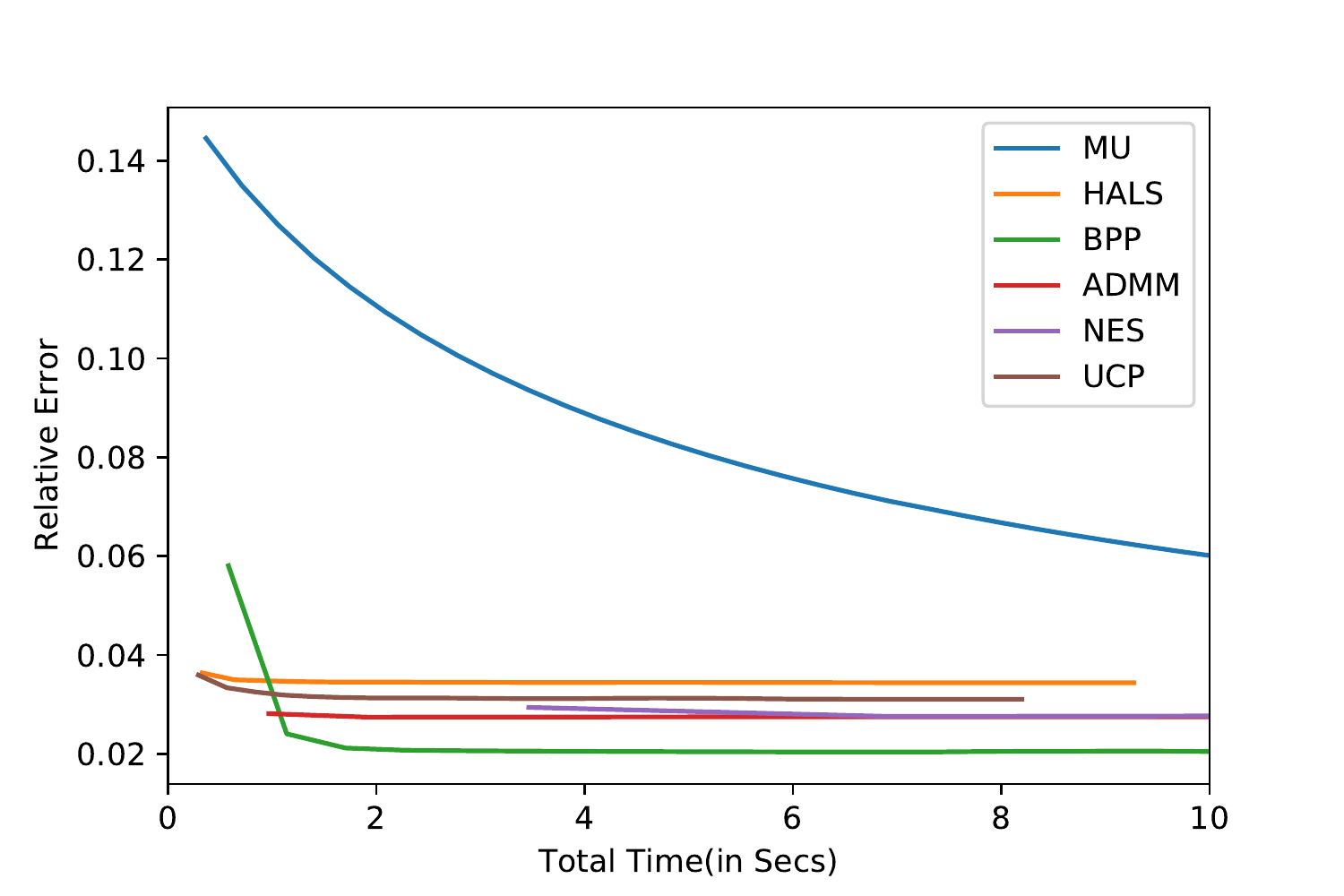}
	}
	\subfloat[Low rank $k=48$ \label{fig:accuracyrwk96}]{
		\includegraphics[width=0.45\textwidth, height=2in]{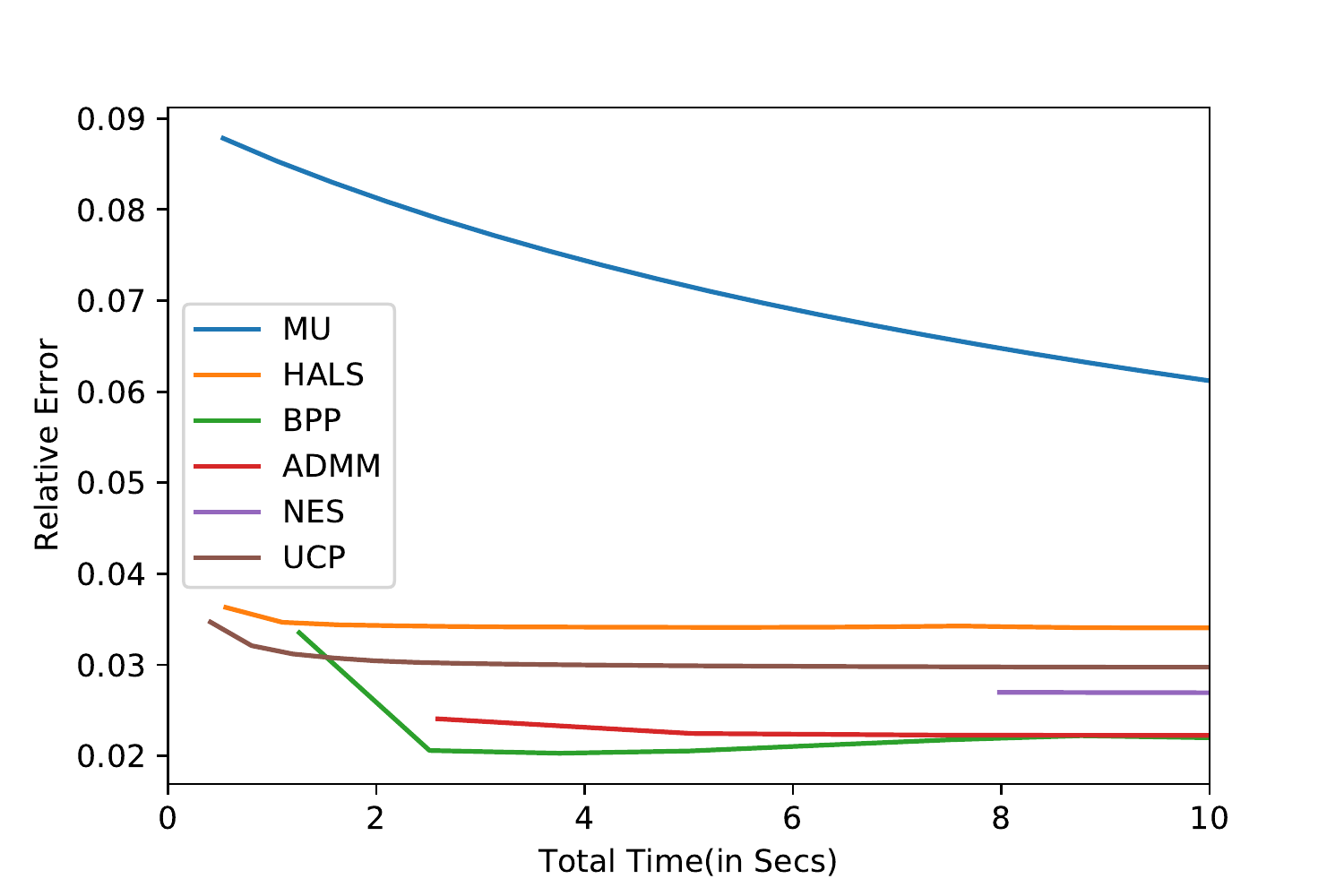}
	}
	\caption{Relative error comparison of updating algorithms on 3D Realworld Low Rank Tensor of size $1447680\times69\times25$ on a 64 Titan Nodes as $64\times1\times1$ Processor Grid for 10 seconds}
	\label{fig:mouseconv}
\end{figure}

\Cref{fig:synconv,fig:mouseconv} show convergence comparisons (error vs time) for each of the updating algorithms on synthetic low-rank and Mouse data sets, using two different target ranks each.
Every algorithm is run for a fixed number of (30) outer iterations for fair comparison.
For the Mouse data in \Cref{fig:mouseconv}, we show only the first 10 seconds because nearly all algorithms are converging within 30 iterations.
The initialized random factors are the same for all algorithms in both tests, and the synthetic tensor is the same for all algorithms.
In both the synthetic and real world cases BPP achieves the lowest approximation error in the shortest amount of time.
Overall results are as expected, such as MU achieving the worst error and ADMM achieving the second best in all cases.
It is also note-worthy that on the real world data set the best algorithms, ADMM and BPP, achieve relative errors of $\approx 2-3 \%$.


\subsection{Scaling Studies}

\subsubsection{Weak Scaling (Synthetic Data)}

\begin{figure}
\centering
\subfloat[Synthetic 3D Low Rank \label{fig:wksca3d}]{
\includegraphics[width=.95\textwidth]{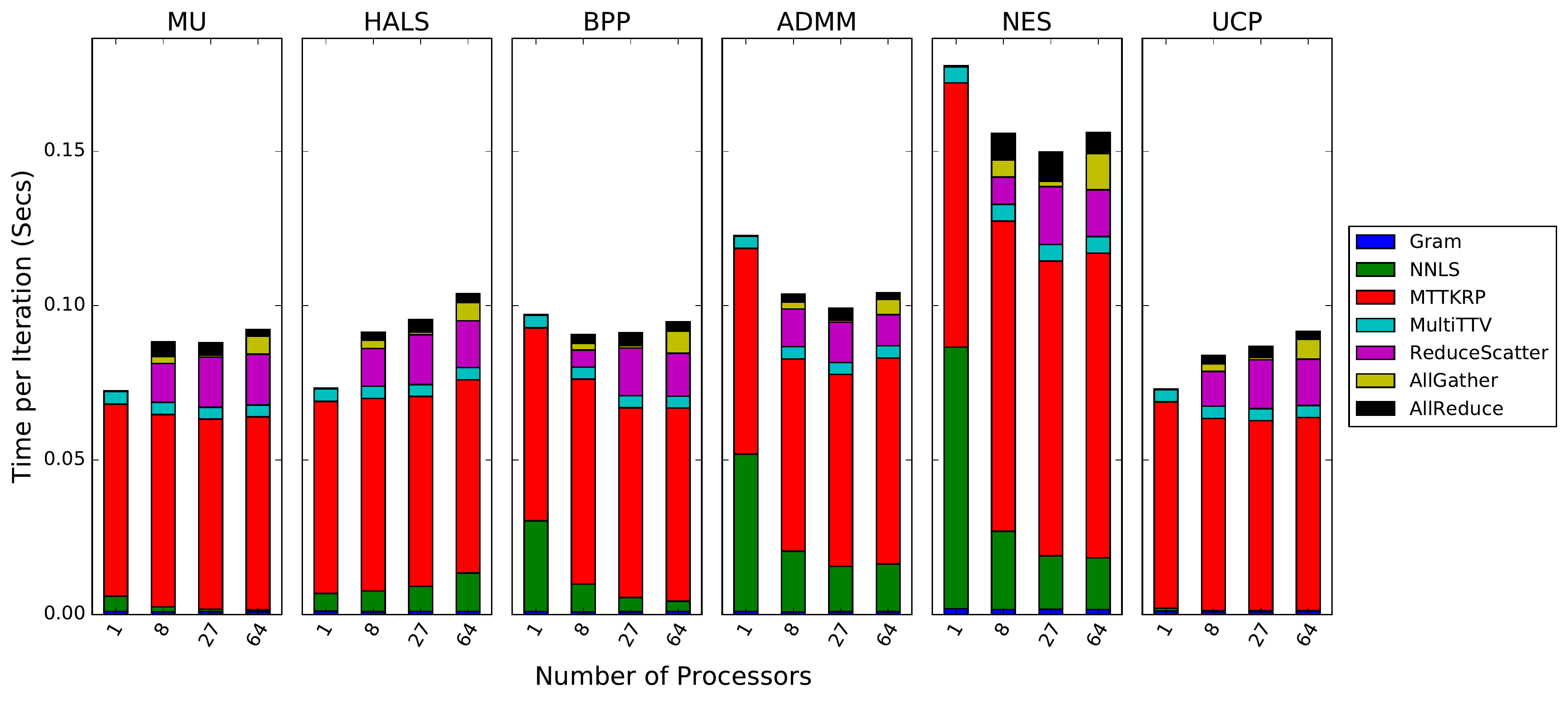}
}\\
\subfloat[Synthetic 4D Low Rank \label{fig:wksca4d}]{
\includegraphics[width=.95\textwidth]{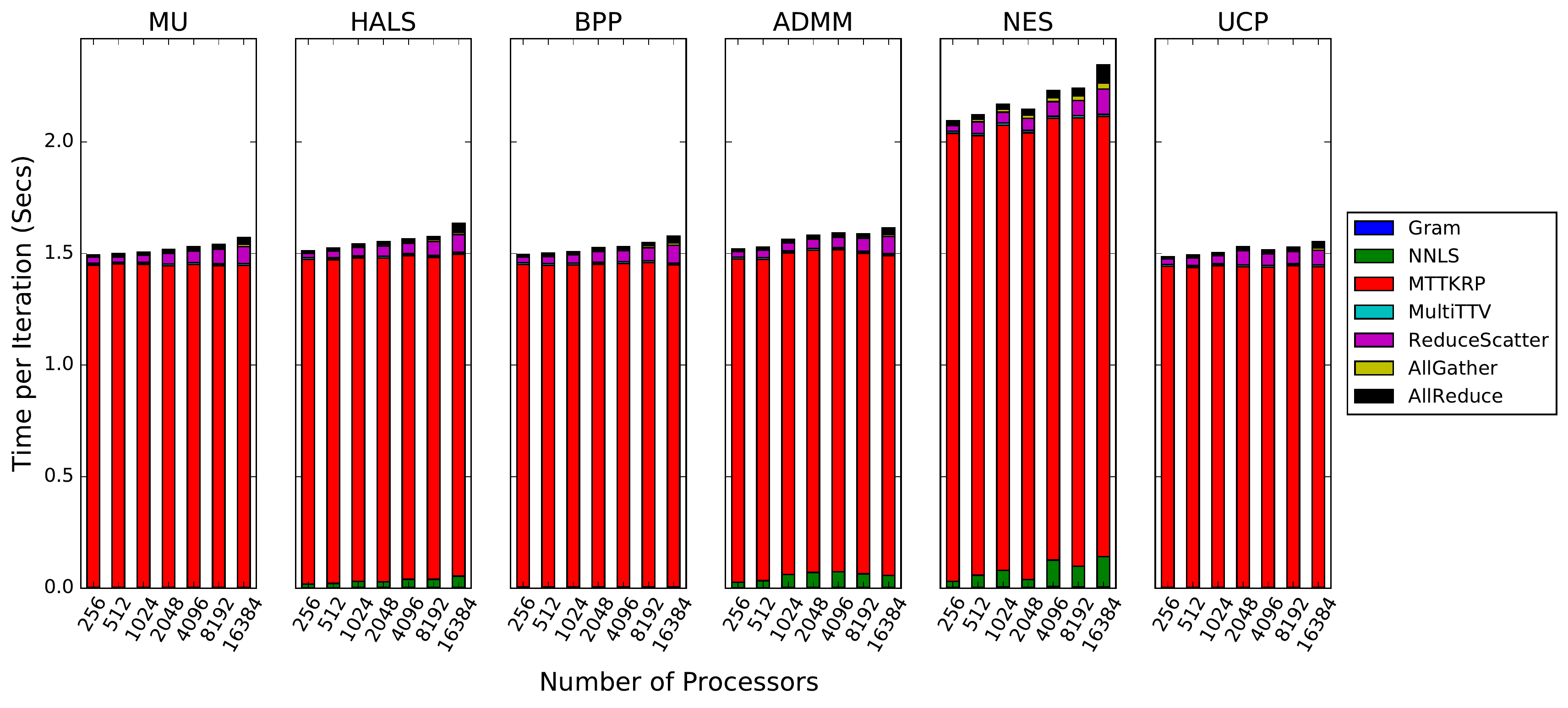}
}
\caption {Weak scaling on synthetic 3D and 4D low-rank tensors. 
For the 3D case, the input tensors are of size $128{\times}128{\times}128$, $256{\times}256{\times}256$, $378{\times}378{\times}378$ and $512{\times}512{\times}512$ on 1, 8, 27 and 64 Titan Nodes. The 4D input tensors are $128{\times}128{\times}128{\times}128$, $256{\times}256{\times}256{\times}256$, $512{\times}512{\times}512{\times}512$, $1024{\times}512{\times}512{\times}512$, $1024{\times}512{\times}1024{\times}512$, $1024{\times}1024{\times}1024{\times}512$, $1024{\times}1024{\times}1024{\times}1024$, $2048{\times}1024{\times}1024{\times}1024$, $2048{\times}1024{\times}2048{\times}1024$
on 1, 16, 256, 512, 1024, 2048, 4096, 8192, and 16384 Titan nodes. 
For all experiments, the low rank is 96.}
\label{fig:synweakscaling}
\end{figure}

We performed weak scaling analysis on 2 different cubical tensors with 3 and 4 modes.~\Cref{fig:synweakscaling} shows the time breakdown for scaling up to 64 nodes of Titan for the 3D case and 16384 nodes for the 4D tensor. In each experiment the size of the local tensor is kept constant at dimension $128$ in each mode for all the runs. 
As expected, the run time is dominated by the cost to compute the MTTKRP, and the domination is more extreme for higher mode tensors. 
Moreover we see reasonable weak scaling as the figures remain relatively flat over all processor sizes. 

The variations occur mainly due to the NNLS and communication portions of the algorithm. These do matter in general and especially for the 3D case where the MTTKRP cost is often comparable to NNLS times, especially for smaller number of processors. However NNLS times scale really well since they split along processor slices rather than fibers and soon become negligible for large processor grids. The amount of communication per processor remains constant but latency costs increase slowly as we scale up.

\subsubsection{Strong Scaling (Synthetic Data)}

\begin{figure}
\centering
\subfloat[Synthetic 3D Low Rank -  $1024{\times}1024{\times}1024$ tensor \label{fig:strsca3d}]{
\includegraphics[width=.95\textwidth]{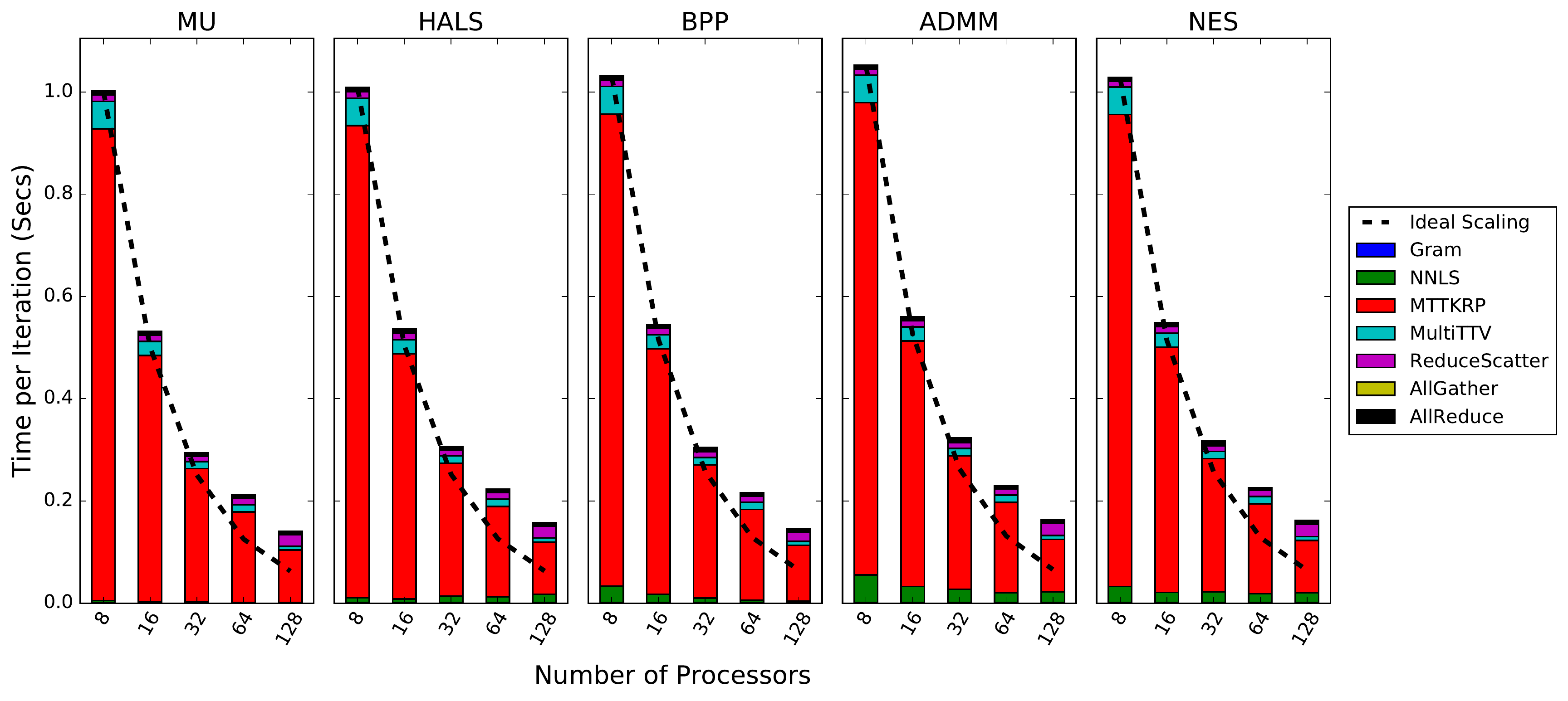}
}\\
\subfloat[Synthetic 4D Low Rank - $512\times512\times512\times512$ tensor \label{fig:strsca4db}]{
\includegraphics[width=.95\textwidth]{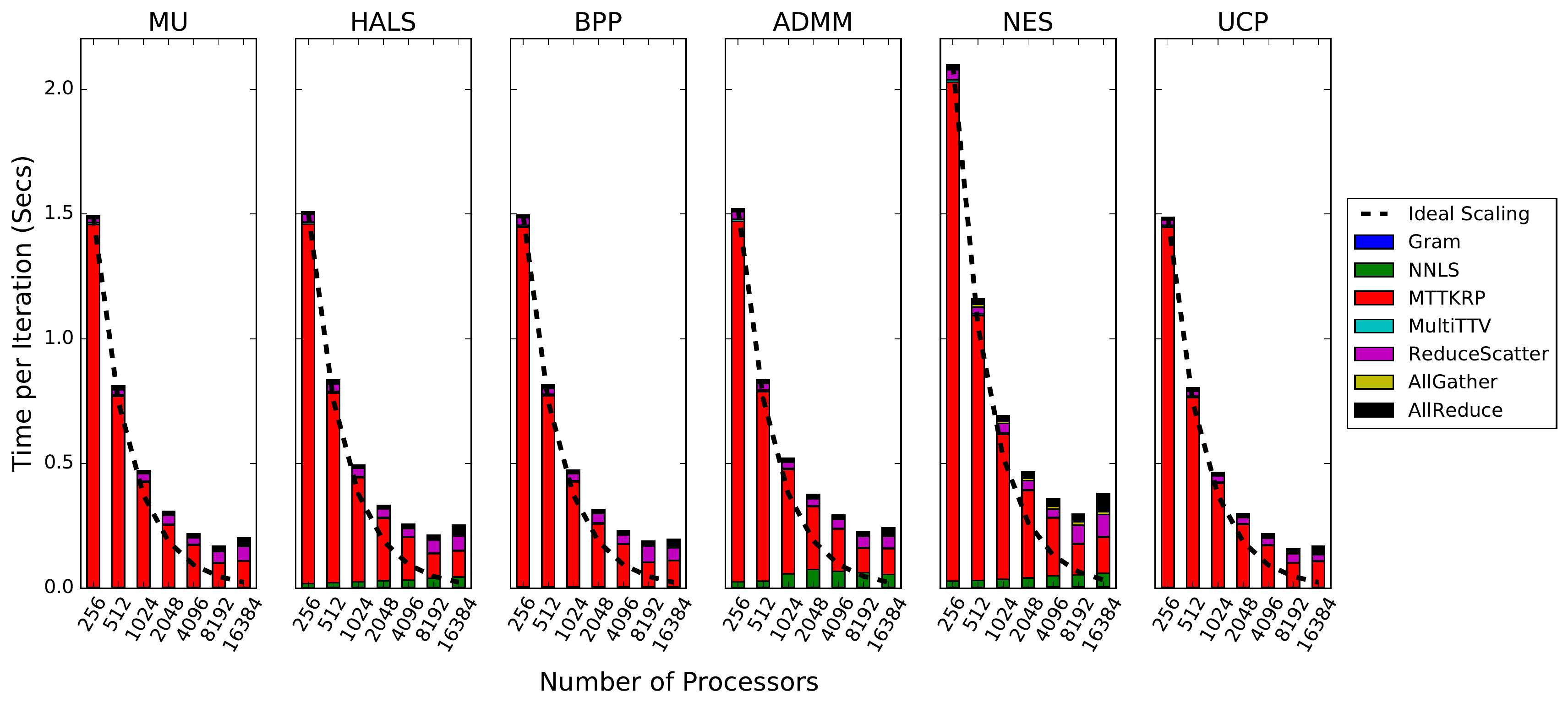}
}
\caption {Strong scaling on synthetic 3D and 4D low rank tensors with low rank 96}
\label{fig:synstrongscaling}
\end{figure}

We run strong scaling experiments on two synthetic square tensors, one 3D and one 4D.~\Cref{fig:synstrongscaling} contains these results for each of the local update algorithms ranging from 1 to 16384 processors. 
We see similar behavior for the 3D and 4D case. For the 3D tensor (\Cref{fig:strsca3d}), we observe good strong scaling up to about 32 nodes and continue to see speed up through 128 nodes. Similarly for the 4D case (\Cref{fig:strsca4db}), the algorithms scale well up to about 1024 nodes and continue to reduce time until 8192 nodes; we observe a slowdown when scaling to 16384 nodes.

One reason for the limit of strong scaling is the communication overheads of AllGather, AllReduce, and ReduceScatter, which become more significant for more processors.
Another reason is that for a cubical tensor of odd dimension, in our case three, the dimension tree optimization is often forced to cast partial MTTKRP into a very rectangular matrix multiplication, depending on the processor grid. Two of the local dimensions must be grouped together while the other is left alone. This means that the largest dimension would need to be close to the product of the other two in order for there to be an approximately square matrix multiplication.



\subsubsection{Strong Scaling (Real World)}

\begin{figure}
\centering
	\includegraphics[width=.95\textwidth]{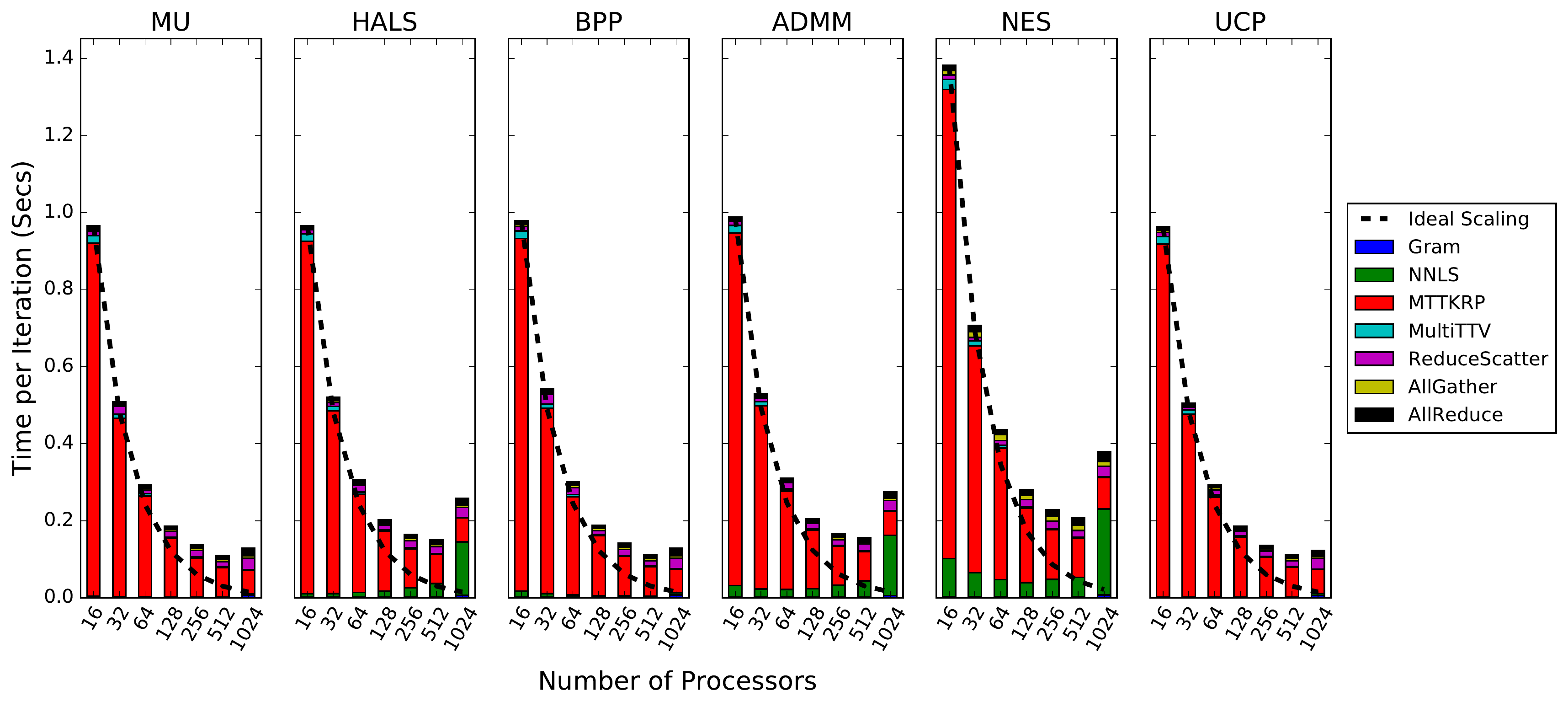}
	\caption {Strong scaling on Mouse dataset}
	\label{fig:mousestrongscaling}
\end{figure}


~\Cref{fig:mousestrongscaling} show strong scaling results on the Mouse dataset. We use a 1D $P\times 1\times 1$ processor grid throughout the experiment.  The results are in line with the synthetic results. We achieve near-perfect scaling up to $\sim32$ nodes and still improve runtimes through 512 nodes. At 1024 nodes the NNLS algorithms, which communicate during the solve steps, perform far worse and show up to $2\times$ slowdown. The non-communicating solvers also degrade in performance but more gracefully.


\subsection{Mouse Data Results}

\begin{figure}
\centering
	\subfloat[Time and trial factors \label{fig:mousecmptt}]{
		\includegraphics[width=.4\textwidth]{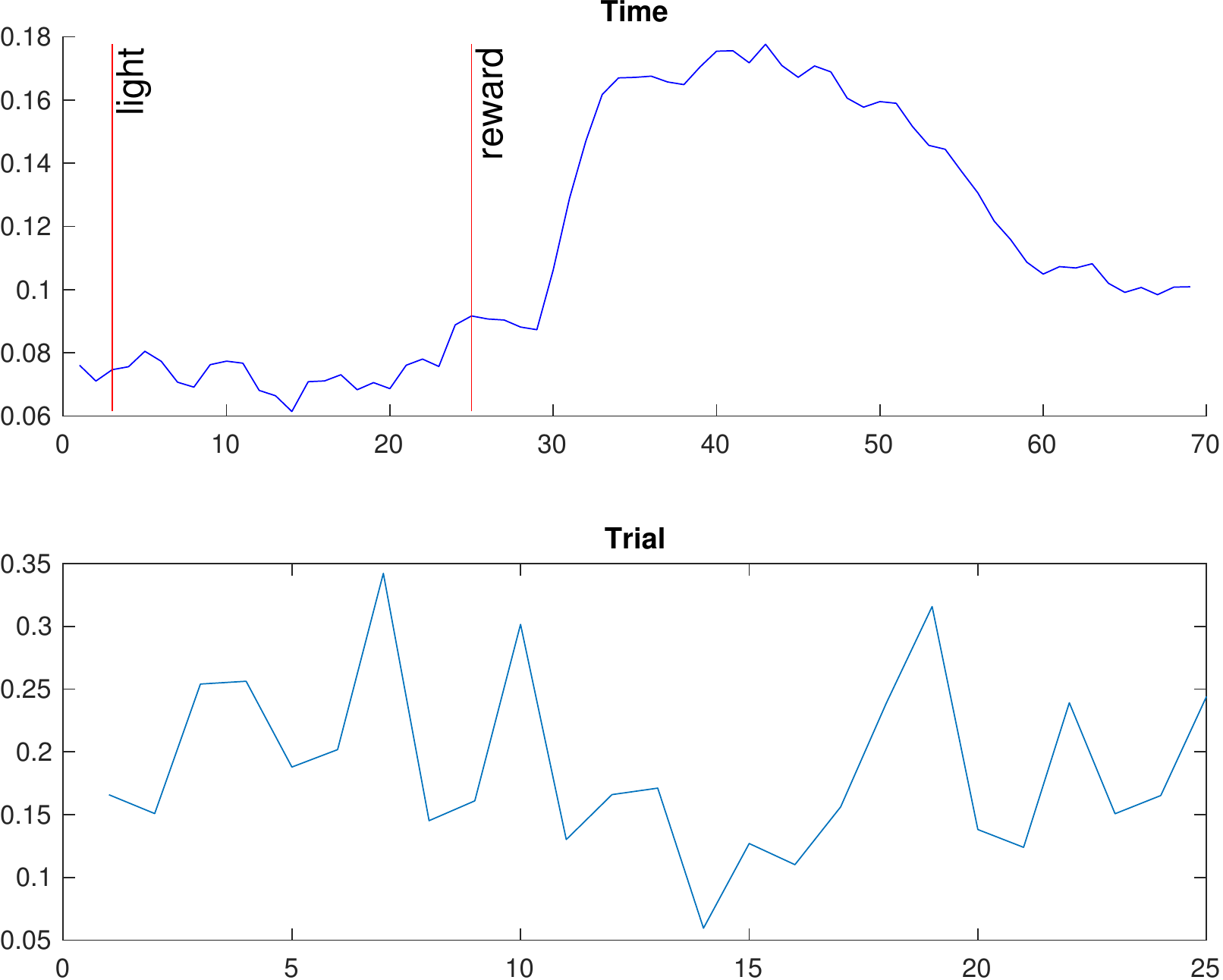}
	}
	\hspace{1cm}
	\subfloat[Brain image factor \label{fig:mousecmpbr}]{
		\includegraphics[width=.4\textwidth]{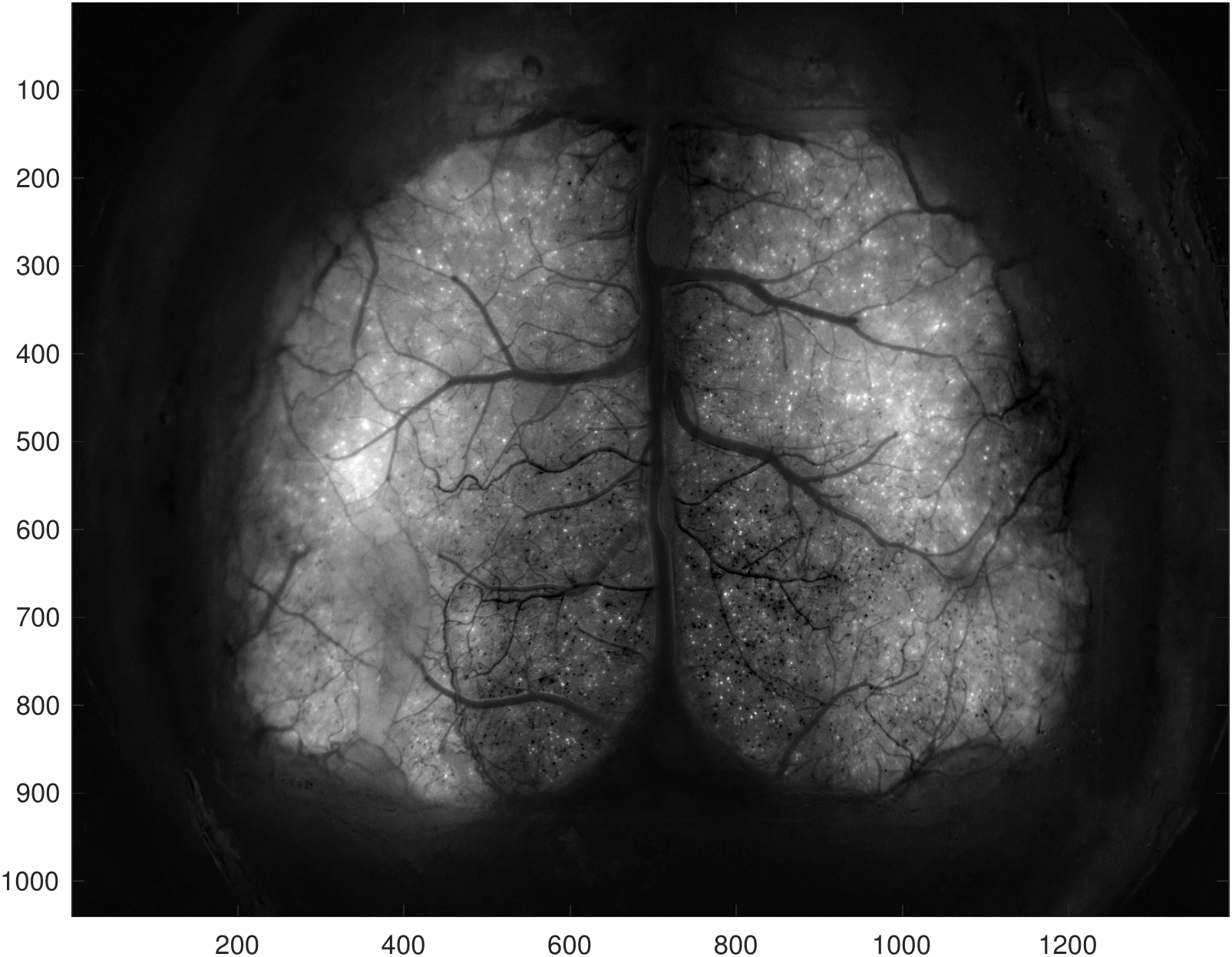}
	}
\caption{Component 22 of rank-32 CP decomposition of Mouse data.}
\label{fig:mousecmp22}
\end{figure}

The CP decomposition of the Mouse data can be used to interpret brain patterns in response to the light stimulus and water reward given to the mouse.
For example, \Cref{fig:mousecmp22} shows a visualization of the factors of the 22nd component of the rank-32 CP decomposition.
From the time factor, we see a marked increase in the importance of the component after the reward time frame, which suggests the activity is a response to the reward.
Because the same mouse undergoes 25 identical trials, we expect to see no pattern in the time factor of each component.
We note that the factors have been normalized, and the absolute magnitude of the y-axis reflects this.
The pixel factor has been reshaped to an image of the same dimensions of the original data.
We observe higher intensity values in the somatosensory cortex (middle, left), which is associated with bodily sensation.
This component, possibly representing a sensory response to the water reward, aligns well with the findings of cell-based analysis \cite[Figure 3]{KZ+16}, which also identified neurons in the somatosensory cortex with intensities that peaked quickly after the reward time frame.

The full set of components for the rank-32 CP decomposition are given in \Cref{fig:mousetimetrials,fig:mousebrains} (\Cref{sec:mouseappendix}).
We note that the interpretation of these computed components is useful only for exploratory analysis.
Their scientific validity would need to confirmed with tests of robustness with respect to choice of rank, random starting point, and algorithm.

\section{Conclusion} \label{sec:conclusion}

In this work, we present PLANC, a software library for nonnegative low-rank factorizations that works for tensors of any number of modes and scales to large data sets and high processor counts.
The software framework can be adapted to use any NNLS algorithm within the context of alternating-updating algorithms.
We use a dimension tree optimization to avoid unnecessary recomputation within the bottleneck local MTTKRP computation, and we use an efficient parallelization algorithm that minimizes communication cost.
Our performance results show the ability to (weakly) scale well on synthetic data to over 16000 nodes (35 TB of data), and we show improved performance by strong scaling on a mouse brain imaging data set of size 20 GB on up to 512 nodes.

PLANC is able to offload some of the computation to a GPU, and we show that this can significant improve the overall runtime.
This is possible because in typical cases, the bottleneck computation (MTTKRP) can be cast as a pair of matrix multiplications (GEMMs) each iteration, which benefit from GPU acceleration for sufficiently large dimensions.
These dimensions depend on the (local) tensor size and the rank of the CP decomposition.
Two of the dimensions can be tuned by the processor grid, which determines the local tensor dimensions, and the choice of dimension tree. 
Therefore two of the matrix multiplication dimensions are typically large.
The third dimension is exactly the rank of the decomposition, so it is typically the smallest dimension and the limitation on GPU efficiency.

The PLANC software framework is designed to be extensible to NNLS algorithms, and we demonstrate how to add an algorithm (the Nesterov-based algorithm) to the library.
While previous work argued that overall performance was agnostic to NNLS algorithm choice \cite{KBP16}, these results show that for NNLS algorithms that involve extra communication or significant computation, the per-iteration running time can be noticeably affected.
In these cases, the time to solution depends both on the per-iteration time and the convergence rate (number of iterations).

PLANC is available at \url{https://github.com/ramkikannan/planc}. It provides both shared and distributed memory parallel algorithms for computing dense NMF, sparse NMF, and dense NTF. Sparse NTF is not currently supported by PLANC but there are plans to provide functionality for sparse NTF in the future.

\section{Acknowledgements}
\small

We thank Tony Hyun Kim and Mark Schnitzer for providing the Mouse dataset and help with interpretation of the CP components.

This material is based upon work supported by the National Science Foundation under Grant No. OAC-1642385 and OAC-1642410.
This manuscript has been co-authored by UT-Battelle, LLC under Contract No. DE-AC05-00OR22725 with the U.S. Department of Energy.  This project was partially funded by the Laboratory Director's Research and Development fund. This research used resources of the Oak Ridge Leadership Computing Facility at the Oak Ridge National Laboratory, which is supported by the Office of Science of the U.S. Department of Energy.

This research used resources of the National Energy Research Scientific Computing Center, a DOE Office of Science User Facility supported by the Office of Science of the U.S. Department of Energy under Contract No. DE-AC02-05CH11231.

 
 \normalsize

\vspace{-0.1in}

\bibliographystyle{ACM-Reference-Format}
\bibliography{paper}




\newpage
\appendix

\section{Full Results for Mouse Data}
\label{sec:mouseappendix}

\begin{figure}
\centering
	\subfloat[Time \label{fig:mousetimes}]{
		\includegraphics[scale=.7]{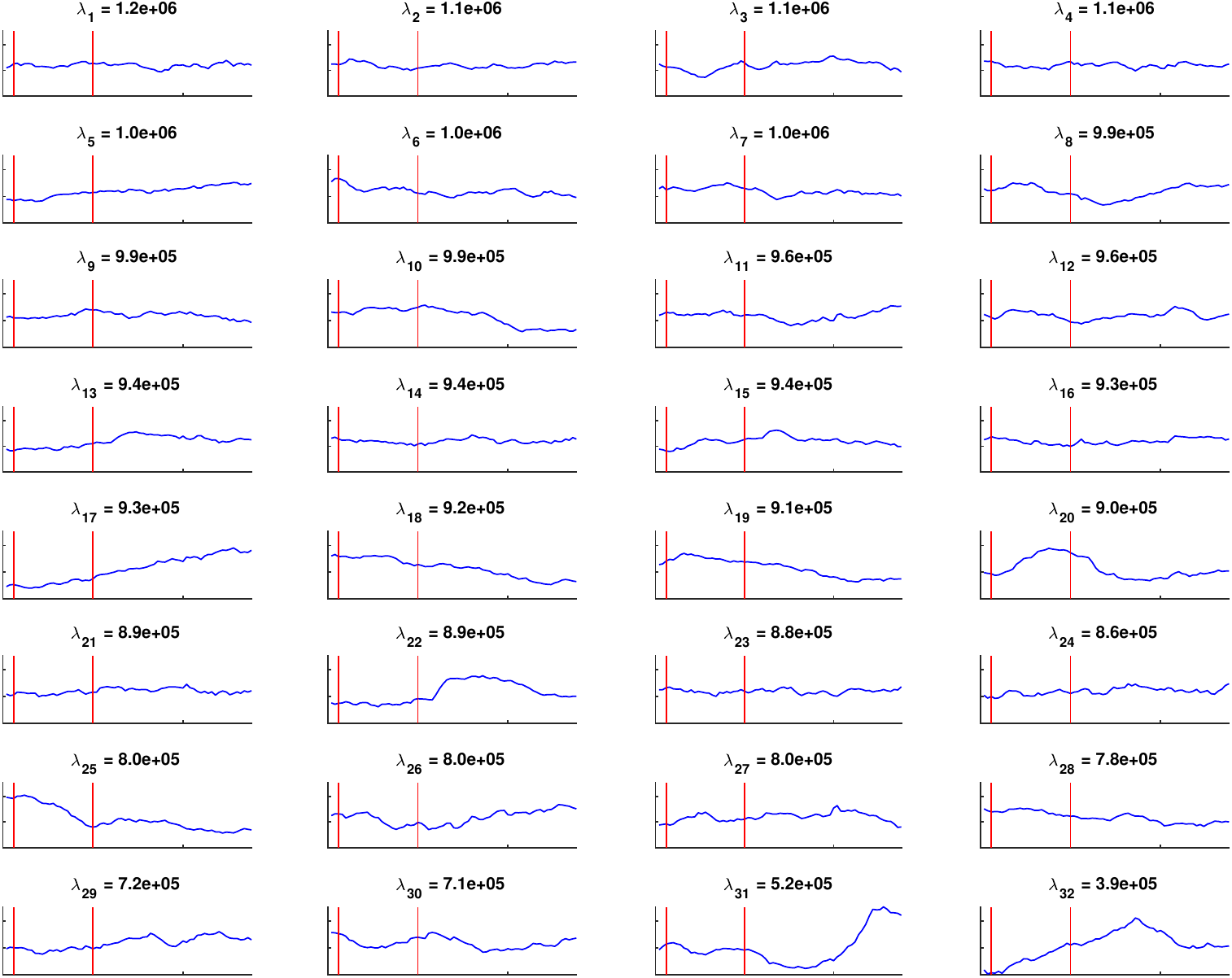}
	} \\
	\subfloat[Trial \label{fig:mousetrials}]{
		\includegraphics[scale=.7]{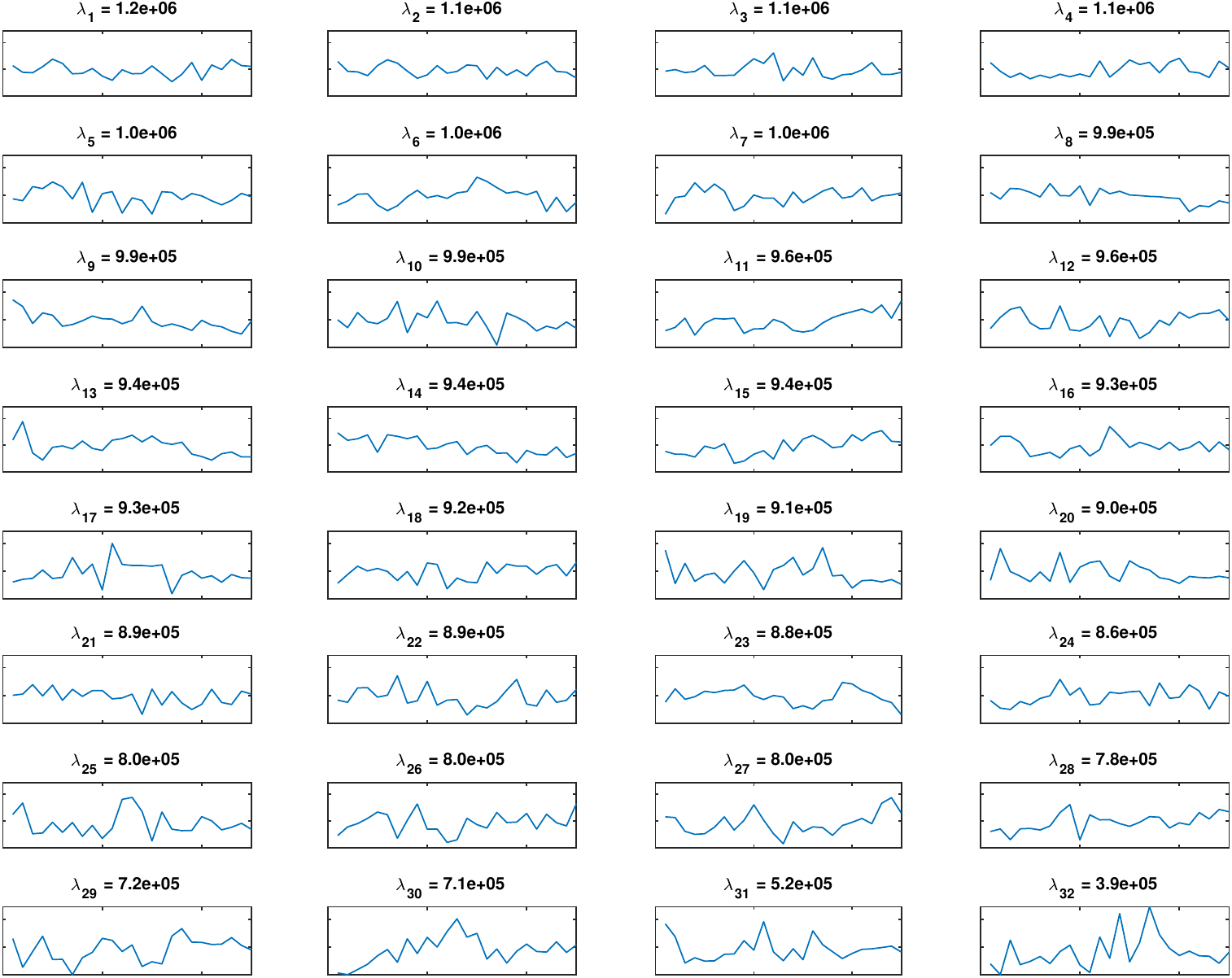}
	}
\caption{Time and trial factors of rank-32 CP decomposition of Mouse data.}
\label{fig:mousetimetrials}
\end{figure}

\begin{figure}
\centering
	\includegraphics[width=.95\textwidth]{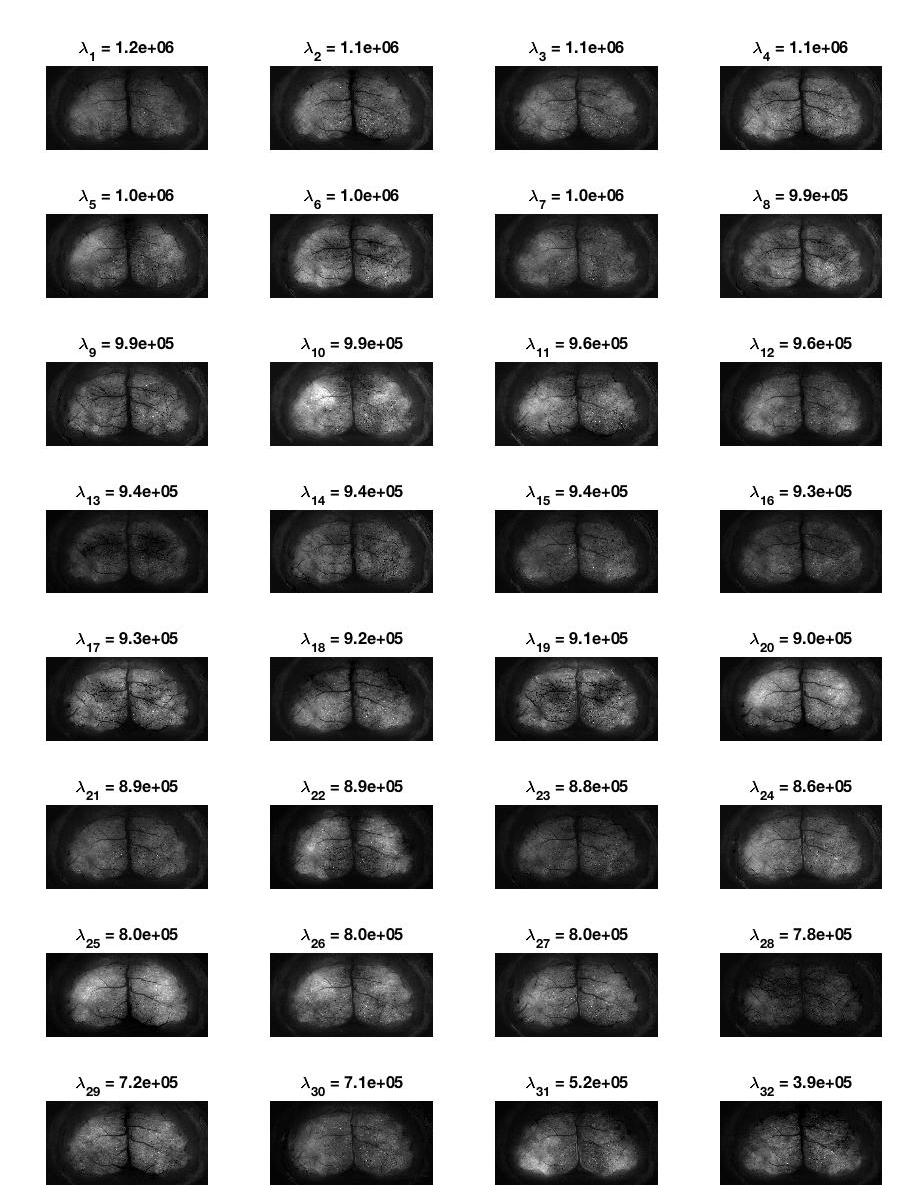}
	\caption{Brain factors of rank-32 CP decomposition of Mouse data.}
	\label{fig:mousebrains}
\end{figure}

\Cref{fig:mousetimetrials,fig:mousebrains} show all 32 components of a CP decomposition of the Mouse data.
This decomposition includes the component highlighted in \Cref{fig:mousecmp22}.
The components are ordered by their weight in the $\V{\lambda}$ vector.
The vertical bars in \Cref{fig:mousetimes} correspond to the frames of the light stimulus and water reward, respectively.

\section{Detailed Parallel Algorithm}

\subsection{Factor Matrix Normalization}
The CP decomposition has a scale indeterminacy.
To prevent possible growth in factor matrix entries, each time a factor matrix is updated, each of the $\Rank$ columns is normalized using the 2-norm and the weights are stored in an auxiliary vector $\V{\lambda}$.
In the distributed algorithm these steps can been seen on lines 21 to 23  in \Cref{alg:Par-NNCP-long}.
Note that communication is required as the global factor matrix norms are computed.

On an algorithmic level one can observe why this step is necessary from the objective function for updating a single factor matrix in the inner iteration $ \min_{\Mn{H}{n}}\| \M{X}_{(n)} - \Mn{H}{n} \M{\Lambda} (\Mn{H}{N} \Khat \cdots \Mn{H}{n+1}  \Khat \Mn{H}{n-1} \cdots \Khat \Mn{H}{1})  \|_{F}^2$, where $\M{\Lambda}$ is the diagonal matrix with the $\V{\lambda}$ vector as its diagonal values.
To solve we simply collapse $\Mn{H}{n} \M{\Lambda}$ together. 
Thus when the solve occurs we are actually computing $\Mn{H}{n} \M{\Lambda}$ which is then normalized to obtain both $\Mn{H}{n}$ and the new $\V{\lambda}$.

\subsection{Relative Error Computation}
\label{sec:error}

Given a model $\T{M}=\CP$, we compute the relative error $\|\T{X} - \T{M}\|/\|\T{X}\|$ efficiently by using the identity $\|\T{X}-\T{M}\|^2 = \|\T{X}\|^2 - 2\langle \T{X}, \T{M} \rangle + \|\T{M}\|^2.$
The quantity $\|\T{X}\|$ is fixed, and the other two terms can be computed cheaply given the temporary matrices computed during the course of the algorithm.
The second term can be computed using the identity $\langle \T{X}, \T{M} \rangle = \langle \Mn{M}{N}, \Mn{H}{N} \rangle$, where $\Mn{M}{N} = \Mz{X}{N} (\Mn{H}{N-1} \Khat \cdots \Khat \Mn{H}{1})$ is the MTTKRP result in the $N$th mode.
The third term can be computed using the identity $\|\T{M}\|^2 = \V{1}^\Tra(\Mn{S}{N} \Hada \MnTra{H}{N} \Mn{H}{N})\V{1}$ where $\Mn{S}{N}=\MnTra{H}{1} \Mn{H}{1} \Hada \cdots \Hada \MnTra{H}{N-1} \Mn{H}{N-1}$.
Both matrices $\Mn{M}{N}$ and $\Mn{S}{N}$ are computed during the course of the algorithm for updating the factor matrix $\Mn{H}{N}$.
The extra computation involved in computing the relative error is negligible.
These identities have been used previously \cite{KB2009,TensorBox,SK16,LKLHS2017}.

\begin{algorithm}
\caption{$(\CPl,\epsilon) = \text{Par-NNCP}(\T{X},R)$}
\label{alg:Par-NNCP-long}
\begin{algorithmic}[1]
\Require $\T{X}$ is an $I_1\times \cdots \times I_N$ tensor distributed across a $P_1\times \cdots \times P_N$ grid of $P$ processors, so that $\T{X}_{\V{p}}$ is $(I_1/P_1)\times \cdots \times (I_N/P_N)$ and is owned by processor $\V{p}=(p_1,\dots,p_N)$, $R$ is rank of approximation
\State \Comment{Initialize data}
\State $a = \text{Norm-Squared}(\T{X}_{\V{p}})$
\State $\alpha = \text{All-Reduce}(a,\textsc{All-Procs})$
\State $\epsilon = $ \texttt{Inf}
\For{$n=2$ to $N$}
	\State Initialize $\Mn{H}{n}_{\V{p}}$ of dimensions $(I_n/P)\times R$ 
	\State $\M[\overline]{G} = \text{Local-SYRK}(\Mn{H}{n}_{\V{p}})$
	\State $\Mn{G}{n} = \text{All-Reduce}(\M[\overline]{G},\textsc{All-Procs})$
	\State $\Mn{H}{n}_{p_n} = \text{All-Gather}(\Mn{H}{n}_{\V{p}},\textsc{Proc-Slice}(n,\VE{p}{n}))$
\EndFor
\State \Comment{Compute NNCP approximation}
\While{$\epsilon > $ \texttt{tol}}
	\State \Comment{Perform outer iteration}
	\For{$n=1$ to $N$}
	\State \Comment{Compute new factor matrix in $n$th mode}
	\State $\M[\overline]{M} = \text{Local-MTTKRP}(\T{X}_{p_1\cdots p_N},\{\Mn{H}{i}_{p_i}\},n)$
	\State $\Mn{M}{n}_{\V{p}} = \text{Reduce-Scatter}(\M[\overline]{M},\textsc{Proc-Slice}(n,\VE{p}{n}))$ 
	\State $\Mn{S}{n} = \Mn{G}{1} \Hada \cdots \Hada \Mn{G}{n-1} \Hada \Mn{G}{n+1} \Hada \cdots \Hada \Mn{G}{N}$
	\State $\Mn[\hat]{H}{n}_{\V{p}} = \text{NNLS-Update}(\Mn{S}{n},\Mn{M}{n}_{\V{p}})$
	\State \Comment{Normalize columns}
	\State $\V[\overline]{\lambda} = \text{Local-Col-Norms}(\Mn[\hat]{H}{n}_{\V{p}})$
	\State $\V{\lambda} = \text{All-Reduce}(\V[\overline]{\lambda},\textsc{All-Procs})$
	\State $\Mn{H}{n}_{\V{p}} = \text{Local-Col-Scale}(\Mn[\hat]{H}{n}_{\V{p}},\V{\lambda})$
	\State \Comment{Organize data for later modes}
	\State $\M[\overline]{G} = {\Mn{H}{n}_{\V{p}}}^\Tra\Mn{H}{n}_{\V{p}}$
	\State $\Mn{G}{n} = \text{All-Reduce}(\M[\overline]{G},\textsc{All-Procs})$
	\State $\Mn{H}{n}_{p_n} = \text{All-Gather}(\Mn{H}{n}_{\V{p}},\textsc{Proc-Slice}(n,\VE{p}{n}))$
	\EndFor
	\State \Comment{Compute relative error $\epsilon$ from mode-$N$ matrices}
	\State $\overline{\beta} = \text{Inner-Product}(\Mn{M}{N}_{\V{p}},\Mn[\hat]{H}{N}_{\V{p}})$
	\State $\beta = \text{All-Reduce}(\overline{\beta},\textsc{All-Procs})$
	\State $\gamma = \V{\lambda}^\Tra (\Mn{S}{N} \Hada \Mn{G}{N}) \V{\lambda}$
	\State $\epsilon = \sqrt{(\alpha-2\beta+\gamma)/\alpha}$ 
\EndWhile
\Ensure $\|\T{X} - \dsquare{\V{\lambda}; \Mn{H}{1},\dots,\Mn{H}{N}}\| /\|\T{X}\| = \epsilon$
\Ensure Local matrices: $\Mn{H}{n}_{\V{p}}$ is $(I_n/P)\times R$ and owned by processor $\V{p}=(p_1,\dots,p_N)$, for $1\leq n \leq N$, $\V{\lambda}$ stored redundantly on every processor
\end{algorithmic}
\end{algorithm}
%

\end{document}